\documentclass[11pt,letterpaper]{amsart}

\usepackage[utf8]{inputenc}
\usepackage[T1]{fontenc}
\usepackage{amsmath, amssymb, amsthm, mathtools}
\usepackage{bbm} 
\usepackage{xcolor}

\usepackage[colorlinks=true, linkcolor=blue, citecolor=blue, urlcolor=blue]{hyperref}

\usepackage[margin=1.2in]{geometry}
\usepackage{tikz}
\usetikzlibrary{arrows.meta}

\DeclareMathOperator{\supp}{supp}

\providecommand{\norm}[1]{\left \lVert#1 \right \rVert}

\newcommand{\CM}{{\mathbb C}}
\newcommand{\NM}{{\mathbb N}}
\newcommand{\RM}{{\mathbb R}}

\newcommand{\TM}{{\mathbb T}}
\newcommand{\ZM}{{\mathbb Z}}

\theoremstyle{plain}
\newtheorem{theorem}{Theorem}[section]
\newtheorem{proposition}[theorem]{Proposition}
\newtheorem{lemma}[theorem]{Lemma}
\newtheorem{corollary}[theorem]{Corollary}

\theoremstyle{definition}
\newtheorem{definition}[theorem]{Definition}

\theoremstyle{remark}
\newtheorem{remark}[theorem]{Remark}

\theoremstyle{definition}

\numberwithin{equation}{section}

\begin{document}
\sloppy
\title[Gaps for monotone potentials]{Large gaps and the ten martini problem for monotone quasiperiodic Schrödinger operators}

\author{Svetlana Jitomirskaya}
\address{Department of Mathematics, University of California,
Berkeley, California 94720}
\email{sjitomi@berkeley.edu}

\author{Tom Stoiber}
\address{Department of Mathematics, University of California,
Berkeley, California 94720}
\address{Simons Institute for the Theory of Computing, Berkeley, CA 94720}
\email{tstoiber@berkeley.edu}

\begin{abstract}
 We prove the ten martini theorem for a large class of monotone potentials. The key strategy is to establish polynomially large gaps for periodic approximations. 
\end{abstract}

\maketitle

\section{Introduction}
	Mark Kac posed the challenge to prove that the almost Mathieu spectrum
is a Cantor set for every irrational frequency. Barry Simon gave it
the name ``ten martini problem'' and popularized it through his
influential lists \cite{simon15,simonXXI}. Its solution \cite{aj}, that relied heavily on the properties of the cosine, naturally led
to the question of how much of this phenomenon itself depends on the
special form of the cosine potential.

We consider quasiperiodic Schr\"odinger operators
\begin{equation}\label{sch}
    (H(x)\psi)_n
    = \psi_{n-1}+\psi_{n+1}+V(x)(n)\psi_n
\end{equation}
on $\ell^2(\ZM)$, where
    $V(x)(n)=f(\{x+n\alpha\}),\;
    f:[0,1)\to\RM.$
The almost Mathieu operator corresponds to
$f(t)=\lambda\cos 2\pi t$, $\lambda\ne0$.
For a given $f$, the corresponding ten martini
problem is to prove that the spectrum is a Cantor set
for \emph{all} irrational $\alpha$\footnote{The spectrum of $H(x)$ is a.e. $x$ independent and, in fact, the same for all $x$ in the context of Theorem~1.1, by minimality and right continuity.}. 
An a.e. statement can miss some of the most interesting features
of the arithmetic phase diagram: transitions in both spectral
type and quantum dynamics that occur on zero measure sets
of frequencies \cite{gj,jliu,jlt,jz}.
 Moreover, $\alpha$ represents magnetic flux
in the corresponding two-dimensional model (with hopping
determined by $f$) \cite{peierls}, and the object of interest
is the entire flux-dependent spectral picture \cite{hofstadter}, the experimental study of which has recently
advanced considerably \cite{nuckolls2025}. This gives
a physics motivation for results valid for every irrational
$\alpha$.  Analytic/smooth $f$ are a natural place to test the
robustness and universality of the ten martini theorem. In this work we study another natural
and popular class: sampling functions $f$ that are monotone on $[0,1)$ (and are, by default, discontinuous, when extended periodically). 

Interest in this class goes back to the Maryland
model \cite{gfp,fp,simonmaryland}, that was developed in the study of quantum kicked rotor, and its early generalizations
\cite{bls}. \footnote{For more recent work see e.g. \cite{jlmaryland,jk1,kach,jymaryland,jk2,jk3,
kps,jls} as well as multidimensional and
other extensions in \cite{cao,csz,kkps, kpsseries,kpsgeneric,kpspert,kz}.}The power of monotonicity  also has a striking dynamical counterpart in
Avila--Krikorian's theory of monotonic cocycles \cite{akmon}.
As discussed in \cite{jk1}, allowing discontinuities makes
a local version of phase monotonicity possible in the
Schr\"odinger case.

We prove the ten martini theorem for a large class of such $f$.
 A simple case of our main result (Theorem~\ref{th:main}) is the following.
\begin{theorem}
\label{th:intro1}
Suppose $f \in C^1([0,1])$ and $f'>0$. Then the  spectrum
of $H$ is a Cantor set for every irrational $\alpha$.
\end{theorem}

    In particular, the theorem applies to the sawtooth potential
$f(t)=\lambda t$ for every $\lambda>0$. The full result, Theorem~\ref{th:main}, requires
considerably less regularity and also allows $f$ that are
unbounded above, such as, for example,
$
f(t)=\frac{\lambda}{1-t}.
$

For the potentials of Theorem~\ref{th:intro1}, the spectrum necessarily has
positive measure~\cite{jk2}. Thus the celebrated zero-measure arguments that
give Cantor spectrum for the critical almost Mathieu operator
and some piecewise constant $f$ \cite{hs,last,ak,dl,bjk}
are unavailable here.

Cantor spectrum for the almost Mathieu operator was originally approached
separately from the Liouville and Diophantine sides.
On the Liouville side, Choi--Elliott--Yui~\cite{ChoiElliottYui}
obtained exponentially small lower bounds on the gaps of rational
approximants. Combined with quantitative spectral continuity,
these gave Cantor spectrum for sufficiently Liouville frequencies.
On the Diophantine side, the approaches used KAM reducibility \cite{eli} and resonant
intersections of eigenvalue branches 
\cite{sinai,gs}, under perturbative assumptions and/or further non-explicit restrictions on $\alpha.$  Puig~\cite{puig} combined duality,
arithmetic localization \cite{j}, and reducibility to obtain Cantor spectrum non-perturbatively for every Diophantine frequency.

Neither Liouville nor any of the above Diophantine approaches could cover the other arithmetic extreme:
the exponential gap bounds were too small, relative to the
spectral continuity estimates, to reach Diophantine frequencies,
while both localization and reducibility have true arithmetic
obstructions and require Diophantine-type conditions.
Bridging the remaining arithmetic region was the main difficulty
overcome in~\cite{aj}. At the time, making the two sides meet
seemed something of a miracle and for a long time, this remained the {\it only}  ten martini
theorem. An important reason for that is that the rational gap estimates of \cite{ChoiElliottYui} were tied to the special
almost Mathieu structure,  and indeed they have not yet been extended to any other $f$ with positive measure spectrum.

Recently, the all-frequency arguments were developed, stretching the Diophantine approaches to all $\alpha$, by employing, correspondingly, almost reducibility \cite{avilaarc1, avilaarc2,ge} instead of reducibility \cite{ayz,lxz}\footnote{These papers prove a stronger, ``dry'' ten martini problem that asks whether
all gaps allowed by gap labeling are open. There has also been
much other recent progress on this question
\cite{bbl,gwx,bf, ayz2, ge-dry,li-dry};
see \cite{cl,band-sofer} for other related developments}, and almost localization instead of localization\footnote{Along with the  intrinsic symplectic structure of the dual center
developed in \cite{gj}.} \cite{gjy24,gjy,gjynew}. 

This established the ten martini theorem
throughout the global open class of analytic Type I potentials, making
the all-frequency statement robust under analytic perturbations\cite{gjy24,gjy}. The Type I class, introduced in \cite{gjy24}, has also become the setting for recent
all-frequency results on the dry ten martini problem
\cite{ge-dry,li-dry,lxz}.   These ideas also
underlie the recent extension to the supercritical\footnote{That is with positive Lyapunov exponent} part
of the spectrum for all real analytic
potentials \cite{gjynew}.
All these approaches rely heavily on duality.

For monotone $f$ there is no useful duality. Moreover, between discontinuities, all finite-volume eigenvalue branches
increase with $x$ and do not intersect, so the gap-opening arguments based on
resonant intersections of opposite-slope branches cannot
apply. The analytic (almost) reducibility  methods
 are likewise unavailable, so none of the existing approaches work. 

It was a surprise in Kachkovskiy--Parnovski--Shterenberg \cite{kps} that, nevertheless, Cantor result for Diophantine $\alpha$ is, in fact, simple, being a corollary of rank-one stability of localization \cite{jk2}. Indeed, modifying $f$ at a single point by decreasing $f(0)$ preserves the conditions of the localization theorem of \cite{jk2}, yet creates a rank-one perturbation 
of $H(0)$. If the spectrum contained an interval, this would have been impossible for generic rank-one perturbations by a theorem of Gordon \cite{gordon} and
del Rio--Makarov--Simon \cite{dms}. The Cantor spectrum in this case is thus a corollary of two substantial
results: the arithmetic localization  of \cite{jk2} and generic sc spectrum of \cite{gordon,dms}. As before, the arithmetic obstructions to localization \cite{jk3} prevent this
approach from giving an all-frequency theorem.
As with Puig's Diophantine result for the almost Mathieu
\cite{puig}, the challenge remained to reach {\it all}
irrational frequencies.

Originally, we planned to follow the idea of \cite{gjy24,gjy}, to try to somehow adapt the almost localization route. Indeed, almost localization type estimates are readily available here \cite{levi}. However, we were not able to find a corresponding argument, and, in desperation, decided to attempt the rational  gaps route \cite{bs,ChoiElliottYui} for Liouville $\alpha$, hoping to obtain lower bounds  that are exponentially small in the denominator as in \cite{ChoiElliottYui} and then miraculously bridge it with the Diophantine result of \cite{kps}.

To our surprise, it turned out that monotonicity and interlacing give a new source of rational
gaps, with a lower bound of order {\it as large as $1/q$} in each fixed energy
interval, thus allowing to stretch this Liouville side approach to {\it all frequencies.}
We now formulate first a special case of this {\it large gap} estimate.

Let $\Sigma$ be the spectrum of $H(x)$ at some irrational $\alpha$, which does not depend on the phase $x$ under our continuity assumption. For a continued-fraction approximant $p/q$ of $\alpha$,  let $\Sigma_q$ denote the union, over all $x$, of the spectra of the corresponding periodic operators with frequency $p/q$.

\begin{theorem}\label{12}
Suppose $f\in C^1([0,1])$, $f'>0$, and $\alpha\notin \mathbb{Q}.$
Then every open interval $J$ meeting $\Sigma$ contains a gap
of $\Sigma_q$ of length at least $c_J/q$, for some $c_J>0$
and all sufficiently large  $q$.
\end{theorem}

To pass from these periodic gaps to Cantor spectrum, we now only need
quantitative spectral continuity estimates adapted to the jump
discontinuity. These give an approximation error smaller than
$1/q$ allowing to complete the proof.
As before, the argument applies under considerably weaker regularity
assumptions and also allows unbounded potentials. We now
state the full result.

Throughout we will assume the standard logarithmic integrability condition
	\begin{equation}
		\label{eq:log-integrability}
		\int_0^1\log(1+|f(x)|)\,dx<\infty
	\end{equation}
that ensures that the Lyapunov exponent is finite.

	\begin{definition}
	We say $f:[0,1)\to\RM$ is $\gamma$-monotone if
	\begin{equation}
		\label{eq:gamma-monotone}
		f(y)-f(x)\geq\gamma(y-x),\qquad 0\leq x<y<1,
	\end{equation}
	with some $\gamma>0$. 
\end{definition}
	
	We also need to make some continuity assumptions. To include unbounded $f$, we express the regularity in terms of the bounded function $(f+i)^{-1}$.
	\begin{definition}
	We say $f:[0,1) \to \RM$ is $\rho$-{\it resolvent-Hölder-continuous}, $0<\rho \leq 1$, if $f$ is continuous and $(f+\imath)^{-1}$ is $\rho$-Hölder-continuous as a function on $(0,1)$. We will denote this by $\rho$-RHC.
	\end{definition}
	If $f$ is resolvent-Hölder-continuous then $f$ is continuous on $[0,1)$ and admits a (possibly infinite) limit at $1.$ 
	\vskip.1in
	Our main result is the ten martini statement for a large class of $f$ satisfying these mild conditions:
    \begin{theorem}\label{th:main}
Suppose $f$ satisfies \eqref{eq:log-integrability}, is bounded from below, $\gamma$-monotone, $\gamma>0$, and is $\rho$-RHC, $\rho>1/2.$ Then $\Sigma$ is a Cantor set for every $\alpha \notin \mathbb{Q}.$ 
\end{theorem}
We also have
\begin{theorem}\label{2}
For $f$ that are $\rho$-RHC, $\rho>0$ and otherwise as in Theoreom~\ref{th:main},  $\Sigma$ is a Cantor set for every non-Diophantine $\alpha$.
\end{theorem}

As discussed above, by the results of \cite{kps},  Cantor spectrum for monotone potentials follows rather easily in the regime of complete localization, using only that localization for monotone potentials is stable under a specific family of rank-one perturbations. In particular, this does not require any Hölder regularity of $f$. Using that argument to cover the Diophantine frequencies, this of course leads to an immediate corollary:
\begin{corollary}\label{c}
 For $f$ that are $\rho$-RHC, $\rho>0$ and otherwise as above, $\Sigma$ is a Cantor set for every $\alpha \notin \mathbb{Q}.$     
\end{corollary}

Even though the H\"older condition in Theorem \ref{th:main} is slightly more restrictive, we present  it the way we do, since it is  interesting even for the  Diophantine case, as it does not require a proof of localization at all and is the first and rather amusing instance of a rational gaps argument stretching to the entire Diophantine regime.
Moreover, the proof of Theorem \ref{th:main} presented here is self-contained modulo the basics, while that of Corollary \ref{c} ultimately depends, through \cite{kps}, on two highly technical results: uniform localization \cite{jk2} and rank-one-generic singular continuous spectrum \cite{gordon,dms}.

\begin{remark}
\begin{enumerate}
\item The assumption that $f$ is bounded from below (or above) is necessary: For all continuous Maryland-type potentials with $f(0+)=-\infty$ and $f(1-)=+\infty$ the spectra $\Sigma_q=\RM$ equal the entire real line for topological reasons \cite{kps}. In fact, as in \cite{kps}, $\Sigma$ is a Cantor set if and only if $f$ is semibounded given the other conditions of the above theorems. \item Resolvent regularity can be relaxed (through the route of Corollary \ref{c}) to at least power-log H\"older, with essentially no change in the argument. 
\item Logarithmic integrability and $\gamma$-monotonicity are standard conditions in the study of monotone potentials.
The latter can be weakened substantially, as in \cite{jls}, see Remark \ref{follow}.
\end{enumerate}
\end{remark}
Both of the above theorems are based on the general form of the large gap estimate
\begin{theorem}
\label{th:gaps}
Suppose $f:[0,1)\to \RM$ is a continuous $\gamma$-monotone function with \eqref{eq:log-integrability},  and $\alpha\notin \mathbb{Q}.$
Then every open interval $J$ meeting $\Sigma$ and such that $L_J=\inf_{E\in J}L(E)>0,$ contains a gap
of $\Sigma_q$ of length at least $c_J/q$, for some $c_J>0$
and all sufficiently large  $q$.

	\end{theorem}
    \begin{remark}
		It is worth contrasting Theorem \ref{th:gaps} with the seminal result of Choi-Elliott-Yui \cite{ChoiElliottYui} who gave exponentially small lower bounds for the almost Mathieu operator. For the noncritical almost Mathieu,
Avila--You--Zhou \cite{ayz} obtained
subexponential lower bounds in $q$ thus stretching the Liouville proof of the dry
ten martini theorem to frequencies admitting any exponentially
accurate rational approximations.
		
		 Our result says nothing about the minimal size of any open gap and in fact, we suspect some gaps may be closed. 
		At the same time, for the {\it largest} gap inside a fixed interval $J$ it is easy to see a posteriori that even the seemingly very good lower bound $c_J/q$ from Theorem~\ref{th:gaps} cannot be optimal: one has Hausdorff convergence of the spectra $\Sigma_q \to \Sigma$   \cite{as,kps}. If $\Sigma$ is a Cantor set then the intersection $\Sigma\cap J$ with any interval has some  gap of  size $\delta>0$, hence $\Sigma_q \cap J$ must also have a gap of size $\delta-o(1)$. The same is true for any quasiperiodic model with Cantor spectrum and qualitative spectral continuity. 

	\end{remark}

We now turn to the quantitative spectral continuity needed
to obtain persistence of gaps for irrational $\alpha$.  Such continuity with respect to frequency goes back to
\cite{ChoiElliottYui}, who obtained a $1/3$-H\"older
estimate. Avron--van Mouche--Simon
\cite{avms} improved the exponent to $1/2$ for general $C^1$
sampling functions.
Better estimates, up to almost matching the modulus of continuity of $f$, became available through (almost) localization, and were obtained in
\cite{jl,JitomirskayaKrasovsky2002,JitomirskayaMavi14}.

For discontinuous $f$ there is of course an obvious obstacle. However, it is manageable. For example, quantitative Hausdorff continuity of the spectrum was obtained in \cite{bbt} for Sturmian operators\footnote{Albeit with distances measured in Farey space.}. More relevant to the present work, in \cite{kps} 
a general result on qualitative Hausdorff continuity was shown for $f$ with a single jump (possibly unbounded). Here we obtain corresponding {\it quantitative}
estimates under the $\rho$-RHC assumption\footnote{This also implies single jump. It is straightforward to replace $\rho$-H\"older by a much weaker modulus of continuity, weakening the conclusion (and thus the set of  $\alpha$ to which the result applies) correspondingly.}. Note that for this we do not require monotonicity.

The proof is a modification of the trial eigenvector strategy as in \cite{avms, JitomirskayaMavi14} and gives the following:
	
	\begin{theorem}
		\label{th:quant_cont}
		Let $f:(0,1)\to R$ be a $\rho$-RHC potential  and $\alpha\notin\mathbb{Q}.$ For any compact  $K\subset \RM$, there exists a constant $C_{\rho,K}$ depending on the $\rho$-Hölder constant, such that for any $E\in \Sigma\cap K$ and all large enough $q$ (depending on $\rho$ and $K$) one has
		$$\mathrm{dist}(E,\Sigma_q) \leq C_{\rho,K}|\Delta_q|^{\frac{\rho}{1+\rho}}.$$
	\end{theorem}
    As usual, we denote here by $q$ a continued fraction denominator of a convergent $p/q$ of $\alpha$ and write
	 $\Delta_q=|\alpha-p/q|.$  	
        We will say $\alpha$ is Diophantine if $\Delta_q>Cq^{-\tau}$ for some $C<\infty,\,\tau>2$. 
	
	As in the continuous case \cite{JitomirskayaKrasovsky2002,JitomirskayaMavi14}, in the regime of positive Lyapunov exponent one can improve the exponent $\frac{\rho}{1+\rho}$ to  almost $\rho.$ 
	\begin{theorem}
		\label{th:lyapunov_controlled_approx} Let $f:(0,1)\to R$ be  $\rho$-RHC potential   and satisfy \eqref{eq:log-integrability}. Let $\alpha$ be Diophantine and $J$ be a compact interval with $L_J =\inf_{E\in J}L(E) > 0$ and such that $E\in J\mapsto L(E)$ is a continuous function. Then there exists for any $0<s<\rho$ a constant $C_{\alpha,s,J}$ such that
		$$\mathrm{dist}(E, \Sigma_q) \leq C_{\alpha,s, J} |\Delta_q|^s$$
		for every $E\in J\cap \Sigma$ and every continued fraction denominator $q$ which is larger than some  $q_0(\alpha,s,J).$ 
		
	\end{theorem}
	\begin{remark} \begin{enumerate}\item Continuity of the Lyapunov exponent is needed for a uniform constant and is automatic in the $\gamma$-monotone setting by the Thouless formula, as the density of states is then absolutely continuous with bounded density \cite{jk2}.
    \item Unlike Theorems~\ref{th:main} and \ref{2} both above Theorems allow $f$ unbounded at $0.$
    \end{enumerate}
    \end{remark}
    
    Both of these spectral continuity results are proven using the same trial eigenvector strategy, the only difference is that positive Lyapunov exponent allows one to construct better trial eigenvectors supported on shorter intervals. The main difficulty is to show that even trial vectors which are supported near a jump discontinuity can be used to bound the distance to the periodic union spectrum, which crucially uses that the sampling function only has a single jump.
	
	Based on those main theorems let us finish the simple proof of the main result:
	\begin{proof}(of Theorem~\ref{th:main})
For monotone potentials the Lyapunov exponent $L(E)$ is positive for Lebesgue-almost every $E$ (see \cite{DamanikKillip2005} and \cite{jk2}).

 Suppose $\Sigma$ contains an interval, then there is, by continuity, some smaller compact interval $J\subset \Sigma$ with $L_J=\inf_{E\in J}L(E)>0$. By Theorem~\ref{th:gaps}, $\Sigma_q$ has some gap of size $c_J/q$ contained in $\operatorname{int}(J)$. Let $E\in\Sigma$ be the midpoint of this gap.
		
		If $\alpha$ is Diophantine and $\rho>\frac{1}{2}$ then there is by Theorem~\ref{th:lyapunov_controlled_approx} for every large enough $q$ some $E'\in \Sigma_q$ with 
		\begin{equation}
			\label{eq:approx_error_dio}
			|E-E'|\leq C_{\alpha, s, J}|\Delta_q|^{s} \leq C_{\alpha, s, J} \frac{1}{q^{2s}}.  
		\end{equation}

		For $\frac{1}{2}<s<\rho$ the right-hand side is eventually smaller than  $\frac{c_J}{2q}$, which is a contradiction with the gap being of size $c_J/q$.
		
		For non-Diophantine $\alpha$ the weaker estimate\begin{equation}
			\label{eq:approx_error}
			|E-E'|\leq C_{\rho, J}\Delta_q^{\frac{\rho}{1+\rho}}.  
		\end{equation} given by Theorem~\ref{th:quant_cont} is already $o(q^{-1})$ along a sufficiently good subsequence of
continued-fraction approximants, for any $\rho>0$.
As above, the midpoint of a gap of $\Sigma_q$ of length
$c_J/q$ inside $J$ gives a contradiction.
\end{proof}
	
\subsection{Main ideas}	
The large gap estimate comes from combining monotonicity with
interlacing. 
	For a fixed  $f$ and  $\alpha,$ let us distinguish the following Hamiltonians, with $\frac{p}{q}$ being a continued-fraction approximant:
	\begin{itemize}
		\item $H(x)$ is the Hamiltonian on $\ZM$ with frequency $\alpha$ and phase $x\in \TM$. 
		\item $H^+(x)$ is the Hamiltonian on the positive half-line $1,2,...$ with frequency $\alpha$ and Dirichlet boundary condition at $0$.   
		\item $H_{q,\theta}(x)$ for $\theta\in \TM$ is the finite-volume Hamiltonian on the chain $0,...,q-1$ at frequency $\frac{p}{q}$ and phase $x\in \TM$ with the twisted-periodic boundary condition $$\psi(-1)=e^{-\imath 2\pi \theta} \psi(q-1), \qquad \psi(q)=e^{\imath 2\pi \theta} \psi(0).$$
		\item By $D_{q}(x)$ we denote the Dirichlet-Hamiltonian with frequency $\frac{p}{q}$ on the chain $1,...,q-1$ (note that this is only $q-1$ sites). By just $D_q$ we refer to $D_q(0)$ as a short-hand.
		\item By $H_{q, \theta}(x+)$ and $H_{q, \theta}(x-)$ we refer to right and left limits. 		
	\end{itemize}	
	
Due to monotonicity the bands making up $\Sigma_q$ are precisely the intervals $$[\inf_{\theta\in \TM}\underline{E}_k(\theta), \sup_{\theta\in \TM}\overline{E}_k(\theta)]$$
where $\underline{E}_k(\theta)$ is the $k$-th eigenvalue of $H_{q,\theta}(0)$ and $\overline{E}_k(\theta)$ is the $k$-th eigenvalue of $H_{q,\theta}(1-)$. Since the Dirichlet restriction $D_q:=D_q(0)$ to sites $1,...,q-1$ is a principal submatrix of both $H_{q,\theta}(0)$ and $H_{q,\theta}(1-)$ simultaneously, one has eigenvalue interlacing: the energy bands of $\Sigma_q$ are separated by the Dirichlet eigenvalues $\mu_1,...,\mu_{q-1}$ of $D_q$. 
To establish gaps it is therefore sufficient to show significant separation of the band bottoms $\underline{E}_{k+1}(\theta)$ from the corresponding interlaced Dirichlet eigenvalues $\mu_k$. Since the Dirichlet Hamiltonian $D_q$ is a principal submatrix of $H_{q,\theta}$ those eigenvalues are related by an exact Schur complement: the eigenvalues $\underline{E}_1(\theta),...,\underline{E}_q(\theta)$ are the zeros of the meromorphic function $$F_\theta(E)= E - f(0)- \sum_{j}\frac{|\langle u_j, b_\theta\rangle|^2}{E-\mu_j}, \qquad b_\theta=\delta_1+e^{\imath 2\pi \theta}\delta_{q-1}$$ which is expressed entirely in terms of the Dirichlet eigenvalues and the boundary values of the Dirichlet eigenfunctions. Using an averaging argument we will derive the existence of zero-free regions of those functions $F_\theta$ and accordingly gaps of size $c_J/q$ in the periodic union spectrum $\Sigma_q$. 

The required analytic input for this argument is formed by three facts about the Dirichlet eigenvalues and eigenfunctions. The first is that by $\gamma$-monotonicity Dirichlet eigenvalues are separated by gaps of at least $\frac{\gamma}{q}$. The second fact is the combination of an exact algebraic identity and a large deviation principle for the Dirichlet determinant which together imply
$$|u_k(1)u_k(q-1)|\leq e^{-(L(\mu_k)+o(1))q}$$
for every Dirichlet eigenfunction. Crucially, this will imply that $F_\theta$ is not very sensitive to the boundary condition $\theta$.  The third fact is that for the Dirichlet Hamiltonians the first site $\delta_1$ is a cyclic vector, which means that its spectral measure must in the limit $q\to \infty$ assign some non-vanishing mass to $J\cap \Sigma$. Indeed, one has $$\liminf_{q\to \infty}\sum_{k: \mu_k\in J}|u_k(1)|^2 >0.$$
	\begin{remark}\label{follow}
		 The analytic input described above follows from $\gamma$-monotonicity of the eigenvalue curves of all periodic rational approximants. By similar arguments as in our recent paper \cite{jls} one can see that those properties are actually satisfied by a significantly larger class of sampling functions than $\gamma$-monotone ones. It is possible to extend our results to that class, but it requires some further technical input, that will be provided in the follow-up paper.
	\end{remark}

{\noindent\small{{\bf Acknowledgements}
This research was supported by NSF DMS-2052899, DMS-2155211, and Simons
896624. The authors would like to thank Netanel Levi for helpful discussions.}}

{\noindent\small{{\bf AI statement}
The authors developed the ideas and proofs and wrote the paper;
they take full responsibility for its contents.
AI assistance in this work was limited in scope.
ChatGPT~5.5 pointed us to an existing algebraic identity relating
eigenvectors and eigenvalues of Jacobi matrices
(see Remark~3.8). ChatGPT~6  helped refine our existing
gap estimate to the bound $c_J/q$ in Theorem~1.9;
the original estimate already sufficed to prove the main theorem.
We also used ChatGPT in preparing figures, proofreading,
and obtaining feedback on exposition.}}

\section{Band spectra of monotone potentials}
\label{sec:bands}
We will assume the conditions of Theorem~\ref{th:gaps} throughout this section.

Let us first characterize the band edges of the periodic spectrum. If $f$ is bounded then $H_{q, \theta}(0-)$ exists and is a bounded operator, as monotone functions admit left and right limits. If $f(0-)=\infty$ then $H_{q,\theta}(0-)$ is a matrix with one diagonal entry equal to infinity. In that case it makes sense to consider the spectrum of $H_{q,\theta}(0-)$ to be defined by taking its smallest $q-1$ eigenvalues to be the eigenvalues $\mu_1,...,\mu_{q-1}$ of the Dirichlet matrix $D_q$ and define its largest eigenvalue to be $\overline{E}_q(\theta) = +\infty$. In that sense, we then have the following:

\begin{lemma}
\label{lemma:monotone_gap}
For each $\theta\in\TM$, denote the ordered eigenvalues of the two endpoint matrices by
\[
\underline E_{1}(\theta)
 \leq\cdots\leq
\underline E_{q}(\theta),
\qquad
\underline E_{j}(\theta)
 :=
 \lambda_j\!\left(H_{q,\theta}(0)\right),
\]
and
\[
\overline E_{1}(\theta)
 \leq\cdots\leq
\overline E_{q}(\theta),
\qquad
\overline E_{j}(\theta)
 :=
 \lambda_j\!\left(
 H_{q,\theta}\!\left(0-\right)
 \right).
\]

Then $\sigma(H_{q,\theta}(x))$ is disjoint from the interval $(\overline E_{j}(\theta), \underline E_{j+1}(\theta))$ for every $x\in \TM$ and $1\leq j\leq q-1$.
\end{lemma}
\begin{proof}
By monotonicity one has
\begin{equation}
\label{eq:monotone}
\underline E_{j}(\theta) \leq E_{j}(\theta,x) \leq  \overline E_{j}(\theta)
\end{equation}
for $x\in (0,\frac{1}{q})$ and $E_{j}(\theta,x)$ the $j$-th eigenvalue of $H_{q,\theta}(x)$. Here we used that the eigenvalue curves are periodic $E_j(\theta, x)= E_j(\theta, x+ \frac{1}{q})$. 

If $f(0-)<\infty$ then from $0-$ to $0$ the potential undergoes a negative jump at the site $0$, which gives a rank one perturbation
\begin{equation}
\label{eq:h_left_right_limit}
H_{q,\theta}(0) = H_{q,\theta}(0-)- (f(0-)-f(0)) |\delta_0\rangle\langle \delta_0|.
\end{equation}
From monotonicity and the Weyl inequalities we conclude
\begin{equation}
\label{eq:evbandedge_inequalities}
\overline E_{j}(\theta) \leq  \underline E_{j+1}(\theta).
\end{equation}
The same inequality holds if $f(0-)=\infty$, since by Cauchy interlacing (see the following Lemma) one then has
$$\overline E_{j}(\theta)  = \mu_j \leq  \underline E_{j+1}(\theta).$$

Combining \eqref{eq:evbandedge_inequalities} and  \eqref{eq:monotone} we conclude that no eigenvalues can be in the interval between $\overline E_{j}(\theta)$ and $\underline E_{j+1}(\theta)$.
\end{proof}

These spectral bands are interlaced with Dirichlet eigenvalues that do not depend on $\theta$:

\begin{lemma}
\label{lemma:monotone_gap_flow}
Let
\[
\mu_{1}<\mu_{2}<\cdots<\mu_{q-1}
\]
be the eigenvalues of the Dirichlet Hamiltonian $D_q=D_q(0)$. Then the possibly empty interval 
\begin{equation}
\label{eq:swept-fixed-theta-interval}
I_{j,\theta}=\left(
 \mu_{j},
 \underline E_{j+1}(\theta)
\right), \qquad 1\leq j\leq q-1
\end{equation}
does not intersect $\bigcup_{x\in \TM} \sigma(H_{q,\theta}(x))$.

\end{lemma}
\begin{proof}
Define the rank one perturbations
\begin{equation}
\label{eq:negative-rank-one-flow}
K_{q,\theta}(t)
 :=
 H_{q,\theta}(0)-t|\delta_{0}\rangle\langle\delta_{0}|,
\qquad t\in \RM
\end{equation}
and let $\kappa_{1}(t,\theta)
 \leq\cdots\leq
\kappa_{q}(t,\theta)
$ be the eigenvalues of $K_{q,\theta}(t)$.

By our choice of labeling, the site $0$ carries the potential
value $f(0)$ in $H_{q,\theta}(0)$. With respect to
\[
\CM^q
 =\CM\delta_{0}\oplus 
 \ell^2(\{1,\ldots,q-1\}),
\]
we have the block decomposition
\begin{equation*}
K_{q,\theta}(t)
 =
 \begin{pmatrix}
 f(0)-t & b^*_\theta\\
 b_\theta & D_q
 \end{pmatrix},
\qquad
b_\theta
 :=
 \delta_1+e^{\imath 2\pi \theta}\delta_{q-1}.
\end{equation*}

Since $D_q$ is a principal minor of each $K_{q,\theta}(t)$ Cauchy interlacing gives
\begin{equation}
\label{eq:left-endpoint-interlacing}
\kappa_{j}(t,\theta)
 \leq
\mu_{j}
 \leq
\kappa_{j+1}(t,\theta),
\qquad 1\leq j\leq q-1
\end{equation}

for every $t \in \RM$. For $t=0$ the matrix $K_{q,\theta}(t)$ is equal to $H_{q,\theta}(0)$, thus we have 
$$\underline E_{j}(\theta)
 \leq
\mu_{j}
 \leq
\underline E_{j+1}(\theta),
\qquad 1\leq j\leq q-1$$
as special cases of \eqref{eq:left-endpoint-interlacing}.
If $f(0-)$ is finite then the inequalities for $t=f(0)-f(0-)$ give
$$\overline E_{j}(\theta)
 \leq
\mu_{j}
 \leq
\overline E_{j+1}(\theta),
\qquad 1\leq j\leq q-1$$
whereas in the case $f(0-)=\infty$ we have by definition
$$\overline E_{j}(\theta)
 =
\mu_{j}
 \leq
\mu_{j+1}=\overline E_{j+1}(\theta),
\qquad 1\leq j\leq q-2$$
together with $\overline{E}_q(\theta)=\infty$ which handles the case $j=q-1$.

As this gives
$$\overline E_{j}(\theta) \leq \mu_j \leq \underline E_{j+1}(\theta)$$
we see that the gap $(\overline E_{j}(\theta), \underline E_{j+1}(\theta))$ from Lemma~\ref{lemma:monotone_gap} already contains the interval $I_{j,\theta}$.
\end{proof}

\begin{figure}[t]
\centering
\begin{tikzpicture}[x=9.2cm,y=1.08cm]

  \def\axisheight{0.25}

  \def\lowerbottomend{0.89}
  \def\lowerbottommid{1.25}
  \def\lowertopend{1.72}
  \def\lowertopmid{2.08}

  \def\upbottomend{3.58}
  \def\upbottommid{3.25}
  \def\uptopend{4.6}
  \def\uptopmid{4.2}

  \def\muheight{2.66}

  \path[fill=black!10]
    (0,\lowerbottomend)
      .. controls (0.18,0.99) and (0.35,1.1) ..
    (0.50,\lowerbottommid)
      .. controls (0.65,1.1) and (0.82,0.99) ..
    (1,\lowerbottomend)
    --
    (1,\lowertopend)
      .. controls (0.82,1.82) and (0.65,2.01) ..
    (0.50,\lowertopmid)
      .. controls (0.35,2.01) and (0.18,1.82) ..
    (0,\lowertopend)
    -- cycle;

  \draw[thick]
    (0,\lowerbottomend)
      .. controls (0.18,0.99) and (0.35,1.1) ..
    (0.50,\lowerbottommid)
      .. controls (0.65,1.1) and (0.82,0.99) ..
    (1,\lowerbottomend);

  \draw[very thick]
    (0,\lowertopend)
      .. controls (0.18,1.82) and (0.35,2.01) ..
    (0.50,\lowertopmid)
      .. controls (0.65,2.01) and (0.82,1.82) ..
    (1,\lowertopend);

  \path[fill=black!16]
    (0,\upbottomend)
      .. controls (0.18,3.50) and (0.35,3.31) ..
    (0.50,\upbottommid)
      .. controls (0.65,3.31) and (0.82,3.50) ..
    (1,\upbottomend)
    --
    (1,\uptopend)
      .. controls (0.82,4.45) and (0.65,4.32) ..
    (0.50,\uptopmid)
      .. controls (0.35,4.32) and (0.18,4.45) ..
    (0,\uptopend)
    -- cycle;

  \draw[very thick]
    (0,\upbottomend)
      .. controls (0.18,3.50) and (0.35,3.31) ..
    (0.50,\upbottommid)
      .. controls (0.65,3.31) and (0.82,3.50) ..
    (1,\upbottomend);

  \draw[thick]
    (0,\uptopend)
      .. controls (0.18,4.45) and (0.35,4.32) ..
    (0.50,\uptopmid)
      .. controls (0.65,4.32) and (0.82,4.45) ..
    (1,\uptopend);

  \draw[very thick,densely dashed]
    (-0.01,\muheight) -- (1.02,\muheight)
    node[right]
    {$\mu_{j}$};

  \node[below]
    at (0.23,0.77)
    {$\underline E_{j}(\theta)$};

  \node[above]
    at (0.23,1.88)
    {$\overline E_{j}(\theta)$};

  \node[below]
    at (0.23,3.39)
    {$\underline E_{j+1}(\theta)$};

  \node[above]
    at (0.23,4.52)
    {$\overline E_{j+1}(\theta)$};

  \draw[->]
    (-0.055,\axisheight) -- (1.075,\axisheight)
    node[right] {$\theta$};

  \draw[->]
    (-0.035,0.18) -- (-0.035,4.98)
    node[above] {$E$};

  \draw
    (0,\axisheight-0.06) -- (0,\axisheight+0.06)
    node[below=4pt] {$0$};

  \draw
    (0.5,\axisheight-0.06) -- (0.5,\axisheight+0.06)
    node[below=4pt] {$\frac12$};

  \draw
    (1,\axisheight-0.06) -- (1,\axisheight+0.06)
    node[below=4pt] {$1$};

  \node[below=15pt]
    at (0,\axisheight)
    {\scriptsize periodic};

  \node[below=15pt]
    at (0.5,\axisheight)
    {\scriptsize antiperiodic};

  \node[below=15pt]
    at (1,\axisheight)
    {\scriptsize periodic};

\end{tikzpicture}

\caption{
Two adjacent bands of $\Sigma_q$ resolved by boundary condition in the $(\theta,E)$-plane. The upper and lower edge of each band are the curves $\overline E_{j}(\theta)$ and $\underline E_{j}(\theta)$ respectively.
}
\label{fig:Dirichlet-interlaced-bands}

\end{figure}
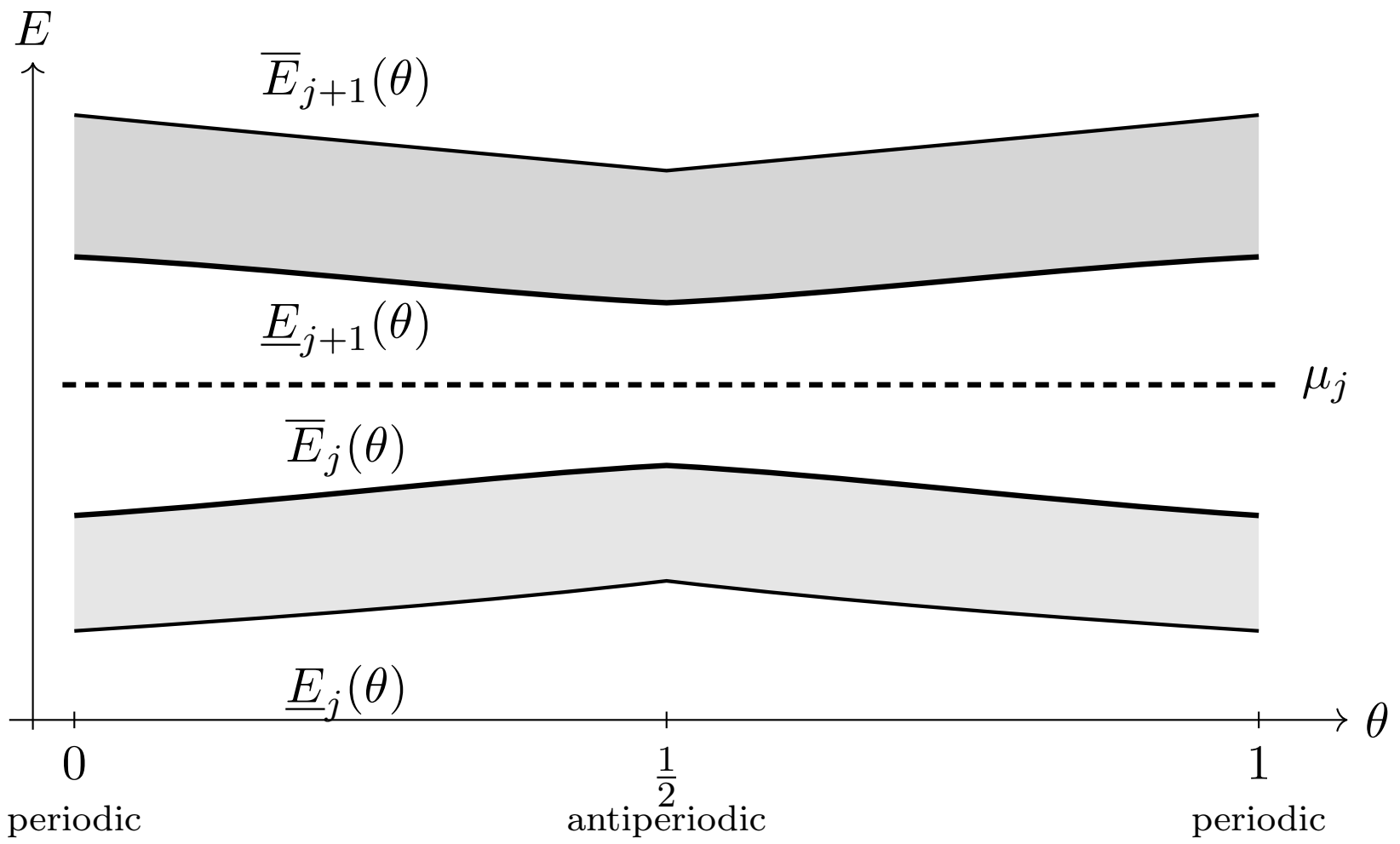

The situation is summarized in Figure~\ref{fig:eigenvalue-curves-dirichlet} and Figure~\ref{fig:Dirichlet-interlaced-bands}.

\begin{figure}[t]
\centering
\begin{tikzpicture}[x=1.25cm,y=1.05cm]
  \def\xmin{-4}
  \def\xmax{4}

  \def\ylowstart{1.05}
  \def\ylowend{2.25}
  \def\yupstart{3.55}
  \def\yupend{4.75}
  \def\ymu{2.75}

  \fill[black!5]
    (\xmin,\ylowend) rectangle (\xmax,\yupstart);
    
  \draw[densely dashed]
    (\xmin,\ymu) -- (\xmax,\ymu)
    node[right=1pt] {$\mu_j$};

  \draw[->]
    (\xmin-0.45,0.55) -- (\xmax+0.55,0.55)
    node[right] {$x$};

  \draw[->]
    (\xmin-0.25,0.45) -- (\xmin-0.25,5.25)
    node[above] {$E$};

  \draw
    (\xmin,0.49) -- (\xmin,0.61)
    node[below=4pt] {$-1/q$};

  \draw
    (0,0.49) -- (0,0.61)
    node[below=4pt] {$0$};

  \draw
    (\xmax,0.49) -- (\xmax,0.61)
    node[below=4pt] {$1/q$};

  \draw[very thick]
    (\xmin,\ylowstart)
      .. controls (-2.75,1.22) and (-1.15,1.82) ..
    (0,\ylowend);

  \draw[very thick]
    (0,\ylowstart)
      .. controls (1.25,1.22) and (2.85,1.82) ..
    (\xmax,\ylowend);

  \draw[very thick]
    (\xmin,\yupstart)
      .. controls (-2.75,3.72) and (-1.15,4.32) ..
    (0,\yupend);

  \draw[very thick]
    (0,\yupstart)
      .. controls (1.25,3.72) and (2.85,4.32) ..
    (\xmax,\yupend);
    
  \filldraw[fill=white]
    (0,\ylowend) circle (1.7pt);

  \filldraw[fill=black]
    (0,\ylowstart) circle (1.7pt);

  \filldraw[fill=white]
    (0,\yupend) circle (1.7pt);

  \filldraw[fill=black]
  (0,\yupstart) circle (1.7pt);

      
  \node[below right]
    at (1.6,3.97)
    {$E_{j+1}(\theta,x)$};

  \node[below right]
    at (-2.7,1.43)
    {$E_j(\theta,x)$};

  \node[right=2pt]
    at (0,\ylowend)
    {$\overline{E}_j(\theta)$};

  \node[left=2pt]
    at (0,\yupstart)
    {$\underline{E}_{j+1}(\theta)$};
\end{tikzpicture}

\caption{
For fixed $\theta$, the eigenvalue branches increase in each periodicity cell and the interval $(\overline{E}_j(\theta), \underline{E}_{j+1}(\theta))$ does not contain any spectrum. The Dirichlet eigenvalues $\mu_j$ interlace the spectrum. In the unbounded case $\mu_j$ coincides with the band edge $\overline{E}_j(\theta)$.
}
\label{fig:eigenvalue-curves-dirichlet}

\end{figure}
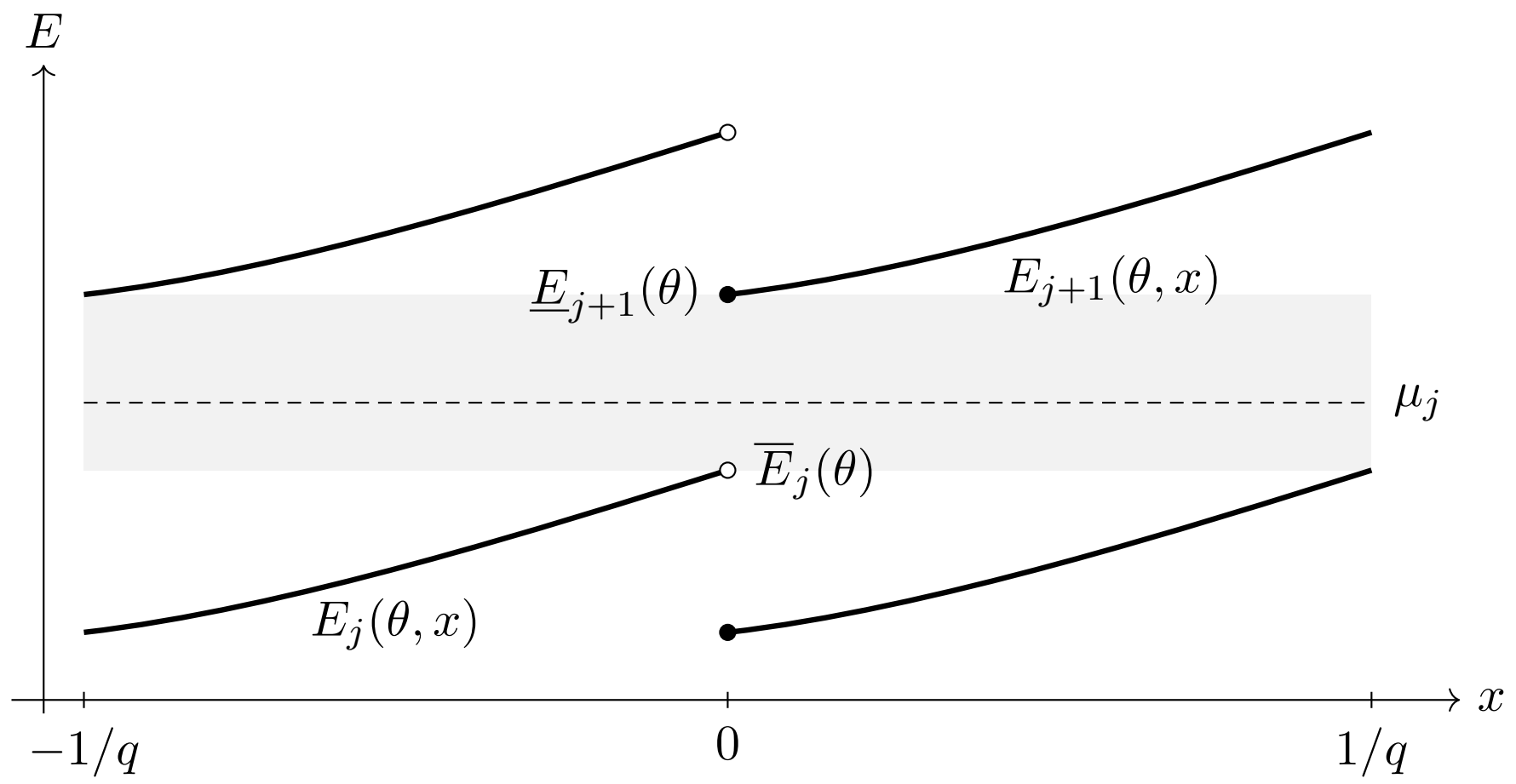

Let us also note immediately a simple but important observation:
\begin{lemma}
\label{lemma:spacing}
We have
\begin{align*}
 \overline{E}_{j+1}(\theta)-\overline{E}_j(\theta) &\geq \frac{\gamma}{q}, \qquad j=1,...,q-1 \\
 \underline{E}_{j+1}(\theta)-\underline{E}_j(\theta) &\geq \frac{\gamma}{q}, \qquad j=1,...,q-1 \\
 \mu_{j+1}-\mu_j &\geq \frac{\gamma}{q}, \qquad j=1,...,q-2
\end{align*}
with $\gamma$ the anti-Lipschitz constant.
\end{lemma}
\begin{proof}
By $\gamma$-monotonicity one has
$$\overline{E}_j(\theta) - \underline{E}_j(\theta) \geq \frac{\gamma}{q}$$
for any $\theta$ since the spectra are $\frac{1}{q}$-periodic. This implies the first two inequalities via
$$\underline{E}_{j}(\theta) \leq \overline{E}_{j}(\theta) \leq \underline{E}_{j+1}(\theta) \leq \overline{E}_{j+1}(\theta).$$
For the Dirichlet eigenvalues one can apply
$$\mu_{j} \leq \underline{E}_{j+1}(\theta) \leq \overline{E}_{j+1}(\theta) \leq \mu_{j+1}.$$
\end{proof}
This statement hinges on the replacement of $\alpha$ by $p/q$; if one simply restricts $H$ with Dirichlet or $\theta$-periodic boundary condition then one does not get similarly strong guaranteed spacing estimates.

\section{Gaps}
\label{sec:gaps}
The strategy to establish spectral gaps is to start with the Dirichlet eigenvalues and try to estimate how far away the $k$-th eigenvalue $\mu_k$ must be from the band bottom $\underline{E}_{k+1}(\theta)$ above it. The distance would be precisely zero if and only if the corresponding Dirichlet eigenfunction $u_k$ satisfies
$$u_k(1) = - e^{\imath 2\pi \theta} u_k(q-1).$$
Indeed, this condition makes the same Dirichlet eigenfunction on $1,...,q-1$ an eigenfunction of $H_{q,\theta}(0)$ on the sites $0,...,q-1$ if one extends it by zero into the site $0$. We can make this quantitative in terms of the endpoint defect which measures how far away any Dirichlet eigenfunction is from satisfying a twisted-periodic boundary condition.

For a normalized eigenvector $u$ of $D_q$ define its endpoint defect by
$$d_q(u):= \left( |u(1)|-|u(q-1)|\right)^2$$

\begin{lemma}
\label{lemma:defect_to_gap}
Fix some compact interval $K\subset \RM$, and let  $\gamma$ be the anti-Lipschitz constant. For all large enough $q$ (depending on $K$) if $u_k$ is a nprmalized eigenvector of the Dirichlet-Hamiltonian $D_q$ with eigenvalue $\mu_k\in K,$ 
then the interval
$$
(\mu_k, \mu_k +\frac{\gamma d_q(u_k)}{5 q}\bigr)
$$
does not intersect $\Sigma_q.$

\end{lemma}
\begin{proof}
We consider the eigenvalue problem of 
\begin{equation*}
H_{q,\theta}(0)
 =
 \begin{pmatrix}
 f(0) & b_\theta^*\\
 b_\theta & D_q
 \end{pmatrix},
\qquad
b_\theta
 :=
 \delta_{1}+e^{\imath 2\pi \theta}\delta_{q-1}.
\end{equation*} By the Schur complement we have
\begin{align*}
\mathrm{det}(E-H_{q,\theta}(0))&=\mathrm{det}(E-D_q) (E- f(0) - b_\theta^* (E-D_q)^{-1} b_\theta)\\
&=\mathrm{det}(E-D_q) \left(E-f(0)- \sum_{j} \frac{|\langle u_j, b_\theta\rangle|^2}{E-\mu_{j}} \right)
\end{align*}
where $u_1,...,u_{q-1}$ are the normalized eigenvectors of $D_q$. Suppose now $E$ is the eigenvalue of $H_{q,\theta}(0)$ closest to $\mu_k$. 

Assume first, for a contradiction, that $E=\mu_k$ and write any eigenfunction of $H_{q,\theta}(0)$ as $\psi = \lambda \delta_0 + \tilde{\psi}$ with $\tilde{\psi}\perp \delta_0$. The eigenvalue equation for $\psi$ is given by the pair of equations
$$\langle b_\theta, \tilde{\psi}\rangle + \lambda (f(0)-\mu_k)=0, \qquad \lambda b_\theta + (D_q-\mu_k)\tilde{\psi}=0.$$

Taking the scalar product of the second with $u_k$ we get $\lambda \langle u_k, b_\theta\rangle=0$. If $\lambda\neq 0$ then this means 
$$0=\langle b_\theta, u_k \rangle= u_k(1) + e^{-\imath 2\pi \theta}u_k(q-1)$$
which implies $d_q(u_k)=0$, a contradiction. If conversely $\lambda=0$ then $\psi$ satisfies Dirichlet boundary conditions at $0$ and since the Dirichlet spectra are simple $\tilde{\psi}$ must be equal to $u_k$ up to a factor. Thus the first equation becomes $\langle b_\theta,u_k\rangle=0$ which again contradicts $d_q(u_k)\neq 0$. 

Next we treat the main case $0<|E-\mu_k|\leq \frac{\gamma}{2q}$. Here we have $\mathrm{det}(E-D_q)\neq 0$ since $E\neq \mu_k$ and the next-closest Dirichlet eigenvalue is at least $\frac{\gamma}{q}$ away from $\mu_k$ by Lemma~\ref{lemma:spacing}. The eigenvalue equation therefore reads
$$0 = E-f(0) - \frac{|\langle u_k, b_\theta\rangle|^2}{E-\mu_{k}}- \sum_{j\neq k} \frac{|\langle u_j, b_\theta\rangle|^2}{E-\mu_{j}}.$$
The next closest eigenvalue $\mu_{j}, j\neq k$ cannot be closer than $\gamma/(2q)$ to $E$, hence
$$\frac{|\langle u_k, b_\theta\rangle|^2}{|E-\mu_{k}|}\leq |E-f(0)| + \frac{2q}{\gamma} \sum_{j\neq k}|\langle u_j, b_\theta\rangle|^2 \leq M_K + \frac{4q}{\gamma}$$
with $M_K=1+\sup_{\mu\in K} |\mu - f(0)|$ and where we used the Bessel inequality with $\|b_\theta\|_2^2=2$. We conclude
$$|E-\mu_{k}|\geq \frac{\gamma |\langle u_k, b_\theta\rangle|^2}{\gamma M_K + 4q} \geq  \frac{\gamma d_q(u_k)}{\gamma M_K + 4q}$$
where the last inequality used that the numerator is minimized when $u_k(1)$ and $e^{\imath 2\pi\theta}u_k(q-1)$ have opposite phase. We then just have to choose $q$ so large that $\gamma M_K<q$.

In the final case $|E-\mu_k|>\frac{\gamma}{2q}$ the same lower bound holds, as one always has 
$$\frac{\gamma}{2q}\geq \frac{\gamma d_q(u_k)}{5q}.$$
In conclusion, for all $\theta$ we have the inclusion
$$\bigl(\mu_k, \mu_k +\frac{\gamma d_q(u_k)}{5 q}\bigr) \subset \bigl(\mu_k, \underline{E}_{k+1}(\theta)\bigr)$$
and the right-hand side is a spectral gap for $H_{q,\theta}$ by Lemma~\ref{lemma:monotone_gap_flow}, hence the interval does not intersect $\Sigma_q$.

\end{proof}

Without further input this gap estimate is not immediately useful, after all, it is to be expected that $u_k(1)$ and $u_k(q-1)$ can both be exponentially small (e.g. in the localization regime when the localization center is close to the center of the interval). More interestingly, it turns out that at positive Lyapunov exponent the product $|u_k(1)u_k(q-1)|$ is always small and hence at least one endpoint weight must be exponentially small:
\begin{lemma}
\label{lemma:endpoint_lyapunov_decay}
For every normalized eigenvalue-eigenfunction pair $(\mu_k, u_k)$ of $D_q$ we have the exact identity
$$|u_k(1)u_k(q-1)|=\frac{1}{|P'_{q}(\mu_k)|} = \prod_{i: i\neq k} \frac{1}{|\mu_i-\mu_k|}$$
with the Dirichlet determinant $P_{q}(E)=\mathrm{det}(D_q-E)$.

Fix any compact interval $J\subset \RM$ with $L_J = \inf_{E\in J}L(E) > 0$. For any $\epsilon>0$ and large enough $q$ we then have
$$|u_k(1)u_k(q-1)|\leq e^{-q (L_J-\epsilon)}$$
for every $\mu_k\in J$.
\end{lemma}
We defer the proof to Section~\ref{sec:proof_lyapunov_lemma} and assume this to hold for now. 

There is now a simple mechanism to explain the existence of polynomial gaps: If $u_k(1)$ is large compared to this exponential Lyapunov scale then its endpoint defect must also be similarly large as $u_k(q-1)$ will be negligibly tiny in comparison. Hence there must be a gap to the right of $\mu_k$ corresponding to this endpoint defect.

The only remaining ingredient is that such eigenfunctions indeed must exist in any interval $J$ for large enough $q$. For this we can use that $\delta_1$ is a cyclic vector for the half-line Schrödinger operator $H^+(0)$ and that the support of its spectral measure contains the bulk spectrum $\Sigma$:

\begin{lemma}
\label{lem:endpoint-spectral-mass}
Let $I\subset\RM$ be an open interval such that $I\cap \Sigma \neq \emptyset$, then there exists a constant $m_I >0$ such that
\begin{equation}
\label{eq:endpoint-mass-liminf}
\liminf_{q\to\infty}
\left\langle
 \delta_1,
 \mathbf 1_I\!\left(D_{q}\right)
 \delta_1
\right\rangle
\geq m_I>0.
\end{equation}
\end{lemma}

\begin{proof}
We can regard $D_q$ as an operator on $\ell^2(\NM)$ by extending it by zero. Note that $D_{q}$ never samples $f$ at the jump discontinuity, hence each diagonal entry of $D_q$ is finite. Likewise, the half-line Dirichlet Hamiltonian at irrational frequency $\alpha$ has finite diagonal entries, but it still is an unbounded operator if $f$ is unbounded. For any $\psi\in \ell^2(\NM)$ which is finitely supported one has the strong convergence $D_{q}\psi \to H^+(0)\psi$ since the potential is pointwise continuous on the finite support of $\psi$. As the vectors of this form are a core of $H^+(0)$, this implies vague convergence of the spectral measures. 

For any continuous function $\chi_I\in C_c(\RM)$ with $\chi_I\leq \mathbf 1_I\!$ we therefore have
$$
\liminf_{q\to\infty}
\left\langle
 \delta_1, \mathbf 1_I\!\left(D_{q}\right)
 \delta_1
\right\rangle
\geq
\lim_{q\to\infty}
\left\langle
 \delta_1, \chi_I\!\left(D_{q}\right)
 \delta_1
\right\rangle
=
\left\langle
 \delta_1,\chi_I(H^+(0))
 \delta_1
\right\rangle=: m_I
$$
and we just need to verify that $m_I>0$ for suitable choice of $\chi_I$. 

Let $\nu^+$ be the spectral measure of $H^+(0)$ at $\delta_1$.
Since $H^+(0)$ is a half-line Jacobi matrix with non-zero off-diagonal coefficients,  $\delta_1$ is cyclic, so the standard Weyl-sequence argument gives
$
\supp\nu^+=\sigma(H^+(0))\supseteq\Sigma.
$
Since $I$ is open and meets $\Sigma$, we can choose
$\chi_I\in C_c(I)$, $0\leq\chi_I\leq1$, such that
$
m_I:=\langle\delta_1,\chi_I(H^+(0))\delta_1\rangle
=\int\chi_I\,d\nu^+>0.
$
\end{proof}

By the mechanism sketched above, the pigeonhole principle already implies the existence of polynomially large gaps. Indeed, the following is a less sharp version of Theorem~\ref{th:gaps} which is nevertheless independently interesting in particular due to its simplicity:

\begin{proposition}
Let $\alpha$ be irrational and $f:[0,1)\to \RM$ be a $\gamma$-monotone increasing continuous function with \eqref{eq:log-integrability}.

Let $J\subset \RM$ be an open interval with $\Sigma\cap J\neq \emptyset$ and such that $L_J=\inf_{E\in J}L(E)>0$. Then there exists a constant $c_J>0$ such that $\Sigma_q$ has a gap of size at least $c_J q^{-2}$ inside $J$ for all large enough continued fraction denominators $q$.    
\end{proposition}
\begin{proof}
Since $J$ is open and $J\cap \Sigma \neq \emptyset$ we can pick a properly contained interval $I\Subset J$ with $I\cap \Sigma \neq \emptyset$. By Lemma~\ref{lem:endpoint-spectral-mass} we have
$$
\sum_{k=1}^{q-1} \mathbf 1_I(\mu_{k})|u_k(1)|^2= \langle \delta_1, \mathbf 1_I(D_q) \delta_1\rangle \geq m_I - o(1)$$
and hence there must be at least one $u_k$ with
$$|u_k(1)|^2\geq \frac{m_I-o(1)}{q}.$$
By Lemma~\ref{lemma:endpoint_lyapunov_decay} applied to $\overline{I}$ one then has
$$d_q(u_k) = |u_k(1)|^2 + |u_k(q-1)|^2 -2 |u_k(1) u_k(q-1)|\geq \frac{m_I-o(1)}{q} - 2e^{-q(L_J-\epsilon)}$$
for all large enough $q$. Lemma~\ref{lemma:defect_to_gap} with $K=\overline{I}$ then turns this into a lower bound for the corresponding band gap. Finally, since $I$ is properly contained in $J$ the gap $\left(\mu_k,\mu_k + \frac{\gamma d_q(u_k)}{5q}\right)$ is contained entirely in $J$ for large enough $q$.
\end{proof}

The gaps of size $q^{-2}$ provided by this simple argument are already enough to prove Cantor spectrum in the non-Diophantine case \footnote{Assuming any resolvent-Hölder regularity}, in fact in case of all $\alpha$ with unbounded continued fraction coefficients. However this is not enough for {\it all} $\alpha$ as one cannot treat the most badly approximated irrationals, indeed, the worst case is the golden ratio with $\Delta_q\sim \frac{1}{\sqrt{5}}q^{-2}$. We therefore aim to obtain a sharper estimate. The best possible gaps one can hope to get with this strategy are of order $q^{-1}$ as that is the ceiling imposed by the eigenvalue spacing. In the following we will prove this estimate.

\subsection{Proof of Theorem~\ref{th:gaps}}

\begin{figure}
    \centering
\begin{tikzpicture}[
    x=1.8cm,y=0.57cm,
    >=Latex,
    branch/.style={black,line width=0.85pt},
    pole/.style={black!55,densely dashed,line width=0.45pt},
    eigenvalue/.style={font=\small,fill=white,inner sep=1.5pt}
  ]
  \pgfmathdeclarefunction{schur}{1}{%
    \pgfmathparse{#1-0.15
      -0.45/(#1+2.0)-0.60/(#1+0.7)
      -0.40/(#1-0.8)-0.55/(#1-2.1)}%
  }

  \foreach \p/\j in {-2.0/1,-0.7/2,0.8/3,2.1/4}{
    \draw[pole] (\p,-4.0)--(\p,4.0);
    \node[below=4pt] at (\p,-4.0) {$\mu_{\j}$};
  }

  \begin{scope}
    \clip (-3.0,-4.0) rectangle (3.08,4.0);
    \foreach \a/\b in {
      -3.0/-2.015,
      -1.985/-0.715,
      -0.685/0.785,
       0.815/2.085,
       2.115/3.08}{
      \draw[branch] plot[domain=\a:\b,samples=240]
        (\x,{schur(\x)});
    }
  \end{scope}

  \draw[->,line width=0.65pt] (-3.12,0)--(3.3,0)
    node[right=3pt] {$E$};
  \draw[->,line width=0.65pt] (-3.12,-4.0)--(-3.12,4.5)
    node[above=3pt] {$F_\theta(E)$};

  \foreach \z/\j in {
      -2.255338042/1,
      -1.125160544/2,
       0.155872622/3,
       1.156166248/4,
       2.418459716/5}{
    \fill (\z,0) circle[radius=1.35pt];
    \node[eigenvalue,below right=4pt] at (\z,0)
      {$\underline E_{\j}(\theta)$};
  }
\end{tikzpicture}
    \caption{Relation of the Schur complement and eigenvalues. If $F_\theta(\mu_k+\delta)$ is negative for some small $\delta$, then $(\mu_k,\mu_k+\delta)$ is contained in a gap of $H_{q,\theta}(0).$}
    \label{fig:schur}
    \end{figure}
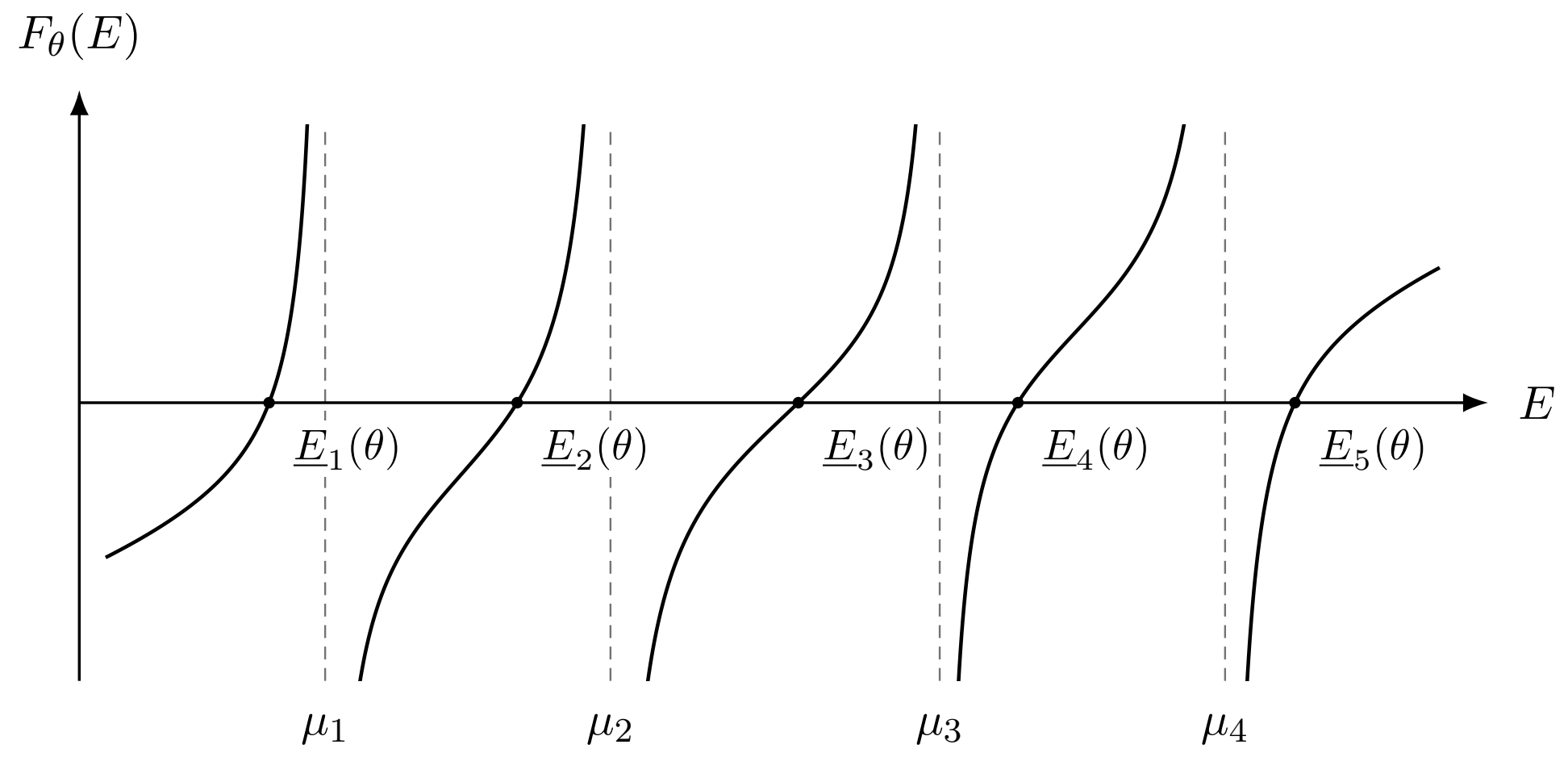

The strategy is similar to Lemma~\ref{lemma:defect_to_gap}: to analyze the eigenvalues of $H_{q,\theta}$ using the Schur complement 
\begin{equation}
\label{eq:Schur_complement}
F_\theta(E)= E - f(0)- \sum_{j}\frac{|\langle u_j, b_\theta\rangle|^2}{E-\mu_j}.
\end{equation}
Generically, its poles are the Dirichlet eigenvalues and its zeros are the eigenvalues of $H_{q,\theta}(0)$, see Figure~\ref{fig:schur}. The strategy for proving gaps is now contained in the following Lemma: 
\begin{lemma}
\label{lem:gapstrategy}
Fix some continued fraction denominator $q$, an open interval $J$ and let $[a,b]\subset J$ be some smaller interval. Suppose there exist positive weights $m_k>0$ for each index $k$ with $\mu_k\in [a,b]$ and some $\delta$ with $0 < \delta < \frac{\gamma}{q}$ such that
$$\sum_{k: \mu_k\in [a,b]} m_k \sup_\theta F_\theta(\mu_k+\delta)< 0,$$
then there exists some index $k$ such that $(\mu_k, \mu_k+\delta)\cap \Sigma_q =\emptyset$.
\end{lemma}
\begin{proof}	
The Schur complement \eqref{eq:Schur_complement} has a pole at each Dirichlet eigenvalue $\mu_j$ except for the (rare) cases in which $\langle u_j, b_\theta\rangle=0$. Note that $F_\theta$ is strictly increasing between poles due to
$$F'_\theta(E)= 1+ \sum_j \frac{|\langle u_j, b_\theta\rangle|^2}{(E-\mu_j)^2}> 0$$
and switches  sign from positive to negative at any pole. Therefore, there is at most one zero on any interval $(\mu_k,\mu_{k+1})$. One can see from the proof of Lemma~\ref{lemma:defect_to_gap} that there are two possible cases for each $j$:
\begin{enumerate}
	\item[(a)] $\underline{E}_{j}(\theta)$ is not a Dirichlet eigenvalue. Then $\underline{E}_{j}(\theta)$ is a zero of $F_\theta$ between $\mu_{j-1}$ and $\mu_j$. Conversely, every zero of $F_\theta$ is an eigenvalue of $H_{q,\theta}(0)$.
	\item[(b)] $\underline{E}_{j}(\theta)$ is a Dirichlet eigenvalue. Then the strict inequalities in $$\mu_{j-2}\leq \underline{E}_{j-1}(\theta)< \overline{E}_{j-1}(\theta) \leq \mu_{j-1} \leq \underline{E}_{j}(\theta) < \overline{E}_{j}(\theta) \leq \mu_j$$ imply that can only have $\underline{E}_{j}(\theta)=\mu_{j-1}$. Equivalently, $\langle u_{j-1}, b_\theta\rangle=0$ which means that $\mu_{j-1}$ is not a pole of $F_\theta$.
\end{enumerate} 
By the assumption there must be some index $k$ such that $F_\theta(\mu_k+\delta)<0$ is negative for every $\theta$. Because of the monotonicity, $F_\theta$ is then negative on the interval $(\mu_k,\mu_k+\delta)$ which implies that this interval is a spectral gap of each $H_{q,\theta}(0)$. 

Let us argue that in fact $\mu_k$ also cannot be an eigenvalue: Suppose for a contradiction that $\underline{E}_{k+1}(\theta)=\mu_{k}$ for some $\theta$, then the given $\theta$ must actually be uniquely determined by the equation $\langle u_{k}, b_\theta\rangle=0$. For any other $\theta'\neq \theta$ the scalar product $\langle u_{k}, b_{\theta'}\rangle$ is non-zero, hence $\mu_k$ is not an eigenvalue of $H_{q,\theta'}(0)$. Therefore, $H_{q,\theta'}(0)$ has no spectrum in $[\mu_k,\mu_k+\delta)$. As one always has $$\mu_k \leq \underline{E}_{k+1}(\theta') < \mu_{k+1}$$ we can conclude $\underline{E}_{k+1}(\theta')>\mu_k+\delta$ for every $\theta\neq\theta'$. But then $\underline{E}_{k+1}(\theta)=\mu_k$ contradicts continuity of the eigenvalues with respect to $\theta'$.

In conclusion, $\underline{E}_{k+1}(\theta)>\mu_k+\delta$ holds for every $\theta$, which implies by Lemma~\ref{lemma:monotone_gap_flow} that $[\mu_k,\mu_k+\delta)$ does not intersect $\Sigma_q$

\end{proof}

If $J$ is an interval with positive Lyapunov exponent $L_J>0$ we will now show that the variation of $F_\theta(E)$ with respect to $\theta$ is negligible due to Lemma~\ref{lemma:endpoint_lyapunov_decay}, so that we can instead study the phase-free version
$$F(E)= E - f(0) - \sum_{j}\frac{m_j}{E-\mu_j}, \qquad m_j = |u_j(1)|^2 + |u_j(q-1)|^2.$$

\begin{lemma}
\label{lem:Ftheta_variation}
Let $[a,b]\subset J$ be properly contained in the bounded open interval $J$ and suppose $L_J>0$. Assume $0<\delta < \frac{\gamma}{2q}$ and set
$$D:=\mathrm{dist}([a,b], \RM\setminus J)$$
then
$$|F(\mu_k+\delta)-F_\theta(\mu_k+\delta)|\leq \frac{8}{D} +(1+\frac{1}{\delta})e^{-(L_J-o(1))q}$$
for all $k$ such that $\mu_k \in [a,b]$ and large enough $q$.
\end{lemma}
\begin{proof}
Assuming without loss of generality that the eigenvectors are real one has 
$$|\langle u_j, b_\theta\rangle|^2 = m_j + 2 \cos(2\pi \theta) u_j(1)u_j(q-1).$$
Lemma~\ref{lemma:endpoint_lyapunov_decay} therefore implies
\begin{align*}
F(\mu_k+\delta)-F_\theta(\mu_k+\delta) &= \sum_{j}\frac{|\langle u_j, b_\theta\rangle|^2-m_j}{\mu_k+\delta-\mu_j}\\
&\leq \sum_{j: \mu_j\notin J}\frac{|\langle u_j, b_\theta\rangle|^2-m_j}{\mu_k+\delta-\mu_j} + 2\sum_{j: \mu_j\in J}\frac{e^{-(L_J-o(1))q}}{|\mu_k+\delta-\mu_j|}
\end{align*}

The second term is exponentially small except for the term with $k=j$ as it is a sum over at most $q$ terms of size $q e^{-(L_J-o(1))q}$ and one can absorb the polynomials into the $o(1)$-loss. 
For large enough $q$ one has $\delta < \frac{1}{2}D$ and therefore $|\mu_k-\mu_j+\delta|> \frac{1}{2}D$ for all $j$ in the first sum, hence
\begin{align*}
|F(\mu_k+\delta)-F_\theta(\mu_k+\delta)| &\leq \frac{2}{D}\sum_j (|\langle u_j, b_\theta\rangle|^2+m_j)+ (1+\frac{1}{\delta}) e^{-(L_J-o(1))q}\\
&=\frac{8}{D} +  (1+\frac{1}{\delta}) e^{-(L_J-o(1))q}
\end{align*}
since the Parseval formula implies $\sum_j |\langle u_j, b_\theta\rangle|^2 \leq 2$ and $\sum_j m_j \leq 2$.
\end{proof}

We therefore only need to prove that a weighted average over $F(\mu_k+\delta)$ is sufficiently negative for suitably small $\delta$. Note that by Lemma~\ref{lem:endpoint-spectral-mass} the sum $\sum_{k: \mu_k \in J} m_k$ admits a (uniformly in $q$) strictly positive lower bound for any interval $J$ with $J\cap \Sigma\neq \emptyset$. This, together with the spacing
$$|\mu_j-\mu_k|\geq \frac{\gamma}{q}|j-k|$$ 
of the Dirichlet eigenvalues will be the only input needed to complete the proof.

\begin{proposition}
\label{prop:F_upper_bound}
Let $J\subset \RM$ be an open interval and fix some intervals $[a,b_0]$, $[a,b_1]$ properly contained in $J$ with $b_0 < b_1$. 

There exists a universal constant $C>0$ such that for every $q$ and every $0<\delta<\frac{\gamma}{2q}$ that is small enough such that
$$\frac{1}{\delta} - C \frac{q^2}{\gamma^2}\delta>0$$
there is some $b \in [b_0,b_1]$ such that 
$$\sum_{k: \mu_k\in [a,b]} m_k F(\mu_k+\delta) \leq  2 |b+\delta-f(0)| - \left(\frac{1}{\delta} - C \frac{q^2}{\gamma^2}\delta\right) \frac{m^2_{[a,b]}}{q} + \frac{8}{b_1-b_0}$$
with $m_{[a,b]}=\sum_{k:\mu_k\in [a,b]} m_k$ 
\end{proposition}
\begin{proof}
We are considering 
$$\sum_{k: \mu_k\in [a,b]} m_k F(\mu_k+\delta)= \sum_{k: \mu_k\in [a,b]} m_k (\mu_k+\delta-f(0))- \sum_{k: \mu_k\in [a,b]} \sum_j \frac{m_k m_j}{\mu_k-\mu_j+\delta}.$$
Since the double sum is almost symmetric in $k$ and $j$ we obtain some cancellations from
$$\frac{1}{\mu_k-\mu_j+\delta}+\frac{1}{\mu_j-\mu_k+\delta}= \frac{-2\delta}{(\mu_j-\mu_k)^2- \delta^2} < 0.$$
In fact, since $|\mu_j-\mu_k|\geq \frac{\gamma}{q}|j-k|$ and $\delta <  \frac{\gamma}{2q}$ we have
\begin{align*}
\sum_{k: \mu_k\in [a,b]} \sum_{j: \mu_j\in [a,b]} \frac{m_k m_j}{\mu_k-\mu_j+\delta} &= \sum_{k: \mu_k\in [a,b]} \frac{m_k^2}{\delta} - \sum_{k<j : \mu_k,\mu_j \in [a,b]} \frac{2\delta m_k m_j}{(\mu_j-\mu_k)^2- \delta^2} \\
&\geq \sum_{k: \mu_k\in [a,b]} \frac{m_k^2}{\delta} - \sum_{k<j : \mu_k,\mu_j \in [a,b]} \frac{2\delta m_k m_j}{\gamma^2 q^{-2} |k-j|^2 - \delta^2}\\
&\geq \sum_{k: \mu_k\in [a,b]} \frac{m_k^2}{\delta} - 4\delta \frac{q^2}{\gamma^2} \sum_{k<j : \mu_k,\mu_j \in [a,b]} \frac{m_k^2 + m_j^2}{|k-j|^2}\\
&\geq \left(\frac{1}{\delta} - C \frac{q^2}{\gamma^2}\delta\right) \sum_{k: \mu_k\in [a,b]} m_k^2
\end{align*}
with some universal constant since
$$\sum_{j\in \ZM\setminus \{k\}}\frac{1}{|k-j|^2} <\infty.$$
The remaining unmatched indices need to be estimated differently. First, we drop the contributions from terms with $k>j$ since those are negative, then we rewrite the sum using indicator functions
$$
-\sum_{k: \mu_k\in [a,b]} \sum_{j: \mu_j\notin [a,b]} \frac{m_k m_j}{\mu_k-\mu_j+\delta} \leq \sum_{a \leq \mu_k \leq b < \mu_j}\frac{ m_k m_j }{\mu_j-\mu_k - \delta}= \sum_{k< j} R(j,k, b)$$
with
$$R(j,k,b)= m_j m_k \frac{1}{\mu_j-\mu_k-\delta} \chi(a\leq \mu_k \leq b) \chi(b < \mu_j).$$
Bounding this naively $R(j,k,b)$ could be as big as $O(q)$ if $\mu_j,\mu_k$ are both close to $b$. However, since we are flexible with the choice of $b$ and allow the upper endpoint to lie somewhere between two real numbers $b_0<b_1$ we can get a better bound by averaging: 
Integrating over the possible upper endpoints we can trivially bound
$$\int_{b_0}^{b_1} R(j,k,b) db \leq m_j m_k \frac{\mu_j - \mu_k}{\mu_j-\mu_k-\delta}\leq 2m_j m_k$$
since the product of the indicator functions is non-vanishing at most for an interval of length $\mu_j-\mu_k$. We then used $|\mu_j-\mu_k-\delta|\geq \frac{1}{2} |\mu_j-\mu_k|$ which is valid for any $\delta<\frac{\gamma}{2q}$.

For the full sum we therefore get
$$\int_{b_0}^{b_1} \sum_{k<j}R(j,k, b) db \leq \sum_{k<j}m_j m_k \frac{\mu_j - \mu_k}{\mu_j-\mu_k-\delta} \leq 2\left(\sum_{k<j} m_k m_j \right) \leq 8.$$
We conclude there must be for any $\delta$ some $b$ such that
$$\sum_{a \leq \mu_k \leq b < \mu_j} m_k m_j \frac{1}{\mu_j-\mu_k - \delta}\leq \frac{8}{b_1-b_0}.$$
If $\delta$ is so small that $\frac{1}{\delta} - C \frac{q^2}{\gamma^2}\delta>0$ then we have for the full expression
\begin{align*}
\sum_{k: \mu_k\in [a,b]} m_k F(\mu_k+\delta)&\leq \sum_{k:\mu_k\in [a,b]} m_k (\mu_k+\delta-f(0)) -  \left(\frac{1}{\delta} - C \frac{q^2}{\gamma^2}\delta\right) \sum_{k: \mu_k\in [a,b]} m_k^2 + \frac{8}{b_1-b_0}\\
&\leq 2 |b+\delta-f(0)| - \left(\frac{1}{\delta} - C \frac{q^2}{\gamma^2}\delta\right) \frac{m^2_{[a,b]}}{q} + \frac{8}{b_1-b_0}
\end{align*}
with $m_{[a,b]}=\sum_{k:\mu_k\in [a,b]} m_k$ where we used that the Cauchy-Schwarz inequality implies for any sum of at most $q$ terms $$q \sum_{k: \mu_k\in [a,b]} m_k^2 \geq  \big(\sum_{k: \mu_k\in [a,b]} m_k\big)^2.$$
\end{proof}

\begin{proof}(of Theorem~\ref{th:gaps})
We have to show that for any interval $J$ with $L_J>0$ and $J\cap \Sigma\neq \emptyset$ there exists some constant $c_J>0$ such that $\Sigma_q$ contains a gap of size $c_J \frac{1}{q}$ within $J$.

By Lemma~\ref{lem:endpoint-spectral-mass} there exists an interval $[a,b_0]$ within $J$ such that $m_{[a,b_0]} > M > 0$ is strictly positive for all large enough $q$. Pick an upper bound $b_1>b_0$ arbitrarily but small enough such that $[a,b_1]$ is still contained in $J$. Set $D=\mathrm{dist}([a,b_1], \RM\setminus J)$.
By Lemma~\ref{lem:Ftheta_variation} the expressions $F(\mu_k+\delta)$ and $F_\theta(\mu_k+\delta)$ satisfy
\begin{align*}
\sum_{k:\mu_k\in [a,b]} m_k\left|F(\mu_k+\delta)- \sup_\theta F_\theta(\mu_k+\delta)\right| &\leq \sum_{k:\mu_k\in [a,b]} m_k \left(\frac{8}{D} + (1+\frac{1}{\delta})e^{-(L_J-o(1))q}\right)\\
&\leq \frac{16}{D} + (2+\frac{2}{\delta})e^{-(L_J-o(1))q}
\end{align*}
uniformly for any choice of $b$ between $b_0$ and $b_1$. By Lemma~\ref{lem:gapstrategy}, to prove we have a spectral gap of the form $(\mu_k, \mu_k+\delta)$ in $\Sigma_q$ we only have to pick $\delta$ so small that the right hand side of
$$\sum_{k:\mu_k\in [a,b]}m_k\sup_\theta F_\theta(\mu_k+\delta) \leq \sum_{k:\mu_k\in [a,b]} m_k F(\mu_k+\delta) + \frac{16}{D} + (2+\frac{2}{\delta})e^{-(L_J-o(1))q}$$
is negative for all large enough $q$. In Proposition~\ref{prop:F_upper_bound} we can make the Ansatz $\delta=\frac{c}{q}$ to obtain
\begin{align*}
\sum_{k:\mu_k\in [a,b]} m_k F(\mu_k+\delta)&\leq  2 |b+\delta-f(0)| - \left(\frac{1}{c} - C \frac{c}{ \gamma^2}\right) M^2 + \frac{8}{b_1-b_0}.
\end{align*}

As required, the expression in brackets is positive for any small enough $c$. Picking $c$ so small that the right-hand side is less than $-\frac{17}{D}$ completes the proof, as after fixing such $c$ the prefactor $2+\frac{2}{\delta}=2+\frac{2q}{c}$ can be absorbed into $e^{-(L_J-o(1))q}$.

\end{proof}

\subsection{Proof of Lemma~\ref{lemma:endpoint_lyapunov_decay}}
\label{sec:proof_lyapunov_lemma}
In this section we assume the setting and notations of Lemma~\ref{lemma:endpoint_lyapunov_decay}. We first prove the algebraic identity
\begin{equation}
\label{eq:endpoint-cofactor-identity}
\left|
 u_{k}(1)u_{k}(q-1)
\right|
= \prod_{i:\, i\neq k} \frac{1}{|\mu_{i}-\mu_{k}|} = \frac{1}{|P_q'(\mu_k)|}
\end{equation}

Let $P_q(E)=\mathrm{det}(D_q-E)$ with normalized eigenvectors $u_1,...,u_{q-1}$, then the adjugate (cofactor) matrix is given by
\begin{align*}
\mathrm{adj}(D_q-E) &= \mathrm{det}(D_q-E) (D_q-E)^{-1}\\
&= \mathrm{det}(D_q-E) \sum_{j} \frac{1}{\mu_{j}-E} |u_j\rangle \langle u_j|  \\
&=\sum_j \prod_{i: \, i\neq j} (\mu_{i}-E) |u_j\rangle \langle u_j|
\end{align*}
for $E$ not in the spectrum since the eigenvalues are all simple. The identity then also holds for all $E$ by continuity. At any eigenvalue $\mu_{k}$ this identity becomes
$$
\operatorname{adj}(D_q-\mu_{k})
 = \left(\prod_{i:\, i\neq k} (\mu_{i}-\mu_k)\right)
 |u_k\rangle\langle u_k| = -P_q'(\mu_k)\,  |u_k\rangle\langle u_k|.
$$
Since $D_q$ is a Jacobi matrix with unit hopping the cofactor which makes up the $(1,q-1)$ entry of the adjugate is, up to sign, the determinant of a triangular matrix with all ones on the diagonal. Hence,
$$
\left|\langle \delta_1, 
 \operatorname{adj}(D_q-E)\delta_{q-1}\rangle
\right|
=1$$
for any $E$. At any eigenvalue $\mu_{k}$ we therefore have the exact identity \eqref{eq:endpoint-cofactor-identity}.

\begin{remark}
\label{rem:eigenvalue_identity}
 The identity is a special case of the off-diagonal eigenvector-eigenvalue identity discussed in \cite[Remark~3]{TaoEtAl}, whose adjugate argument is also what inspired the simple proof given above. However, the identity certainly has a longer history. The earliest statement for tridiagonal matrices seems to be the cofactor formula in \cite{Paige}, and it has also appeared in later works. 
\end{remark}
This product over eigenvalues can now be estimated using the large deviations/Koksma inequality approach of \cite{jk2}. Since we are using Dirichlet determinants at rational frequencies we first need to prove that the eigenvalue counting functions at the single phase $x=0$ approximate the density of states of the irrational frequency bulk Hamiltonian at the required rate of convergence:

\begin{lemma}
\label{lem:rational-irrational-Dirichlet-comparison}
Let $N(E)$ be the integrated density of states for $H$. The eigenvalue counting function 
$$\mathcal N_q^D(E)
 :=
 \mathrm{Tr}\mathbf 1_{(-\infty,E]}
 \bigl(D_q(0)\bigr).
$$
satisfies
\begin{equation}
\label{eq:rational-Dirichlet-counting-discrepancy}
\sup_{E\in\RM}
\left|
 \frac{1}{q-1}\mathcal N_q^D(E)-N(E)
\right|
\leq
\frac{C}{q}
\end{equation}
for a universal constant $C$.
\end{lemma}

\begin{proof}

The result is known already for the Dirichlet restriction $\widehat D_q(x)$ of the
irrational operator $H(x)$ to the chain $\{1,\ldots,q-1\}$
\[
(\widehat D_q(x)\psi)(n)
 =
 \psi(n-1)+\psi(n+1)
 +
 f(x+n\alpha)\psi(n),
 \qquad 1\leq n\leq q-1,
\]
with Dirichlet boundary conditions at $0$ and $q$. We will reduce to this case by comparison of spectra. For the convergent $p/q$ of $\alpha$ set $
\eta_q:=\alpha-\frac pq,$ so that $  \Delta_q=|\eta_q|.$

Define the two real numbers (i.e. without the usual mod $1$)
\[
x_q^{(0)}:=0,
\qquad
x_q^{(1)}:=-(q-1)\eta_q,
\]
and let
\[
x_q^-:=\min\{x_q^{(0)},x_q^{(1)}\},
\qquad
x_q^+:=\max\{x_q^{(0)},x_q^{(1)}\}.
\]
Let us now show that
\begin{equation}
\label{eq:rational-irrational-sandwich}
\widehat D_q(x_q^-)
\leq
D_q(0)
\leq
\widehat D_q(x_q^+).
\end{equation}
For every $1\leq n\leq q-1$ one has
\begin{align*}
x_q^{(0)}+n\alpha
 &=
 n\frac pq+n\eta_q,\\
x_q^{(1)}+n\alpha
 &=
 n\frac pq-(q-1-n)\eta_q.
\end{align*}
It follows that
\begin{equation}
\label{eq:orbit-point-bracketing}
x_q^-+n\alpha
\leq
n\frac pq
\leq
x_q^++n\alpha.
\end{equation}
Note that
$$
\|n\frac pq\|\geq\frac{1}{q}
$$
for $1\leq n\leq q-1$ since $q$ is not a divisor of $np$.
Since
$
|\Delta_q|
<
\frac1{q^2},
$
we obtain both
$$n|\Delta_q| \leq (q-1) |\Delta_q| < \frac{1}{q} $$
and
$$(q-1-n)|\Delta_q| \leq (q-1) |\Delta_q| < \frac{1}{q} $$

Therefore, all three points in
\eqref{eq:orbit-point-bracketing} lie in the same interval
between two consecutive integers, i.e. the arc in between does not cross any jump discontinuity of $f$. Since $f$ is monotone between the jump discontinuities \eqref{eq:orbit-point-bracketing} implies
$$
f(x_q^-+n\alpha)\leq f\left(n\frac pq\right) \leq f(x_q^++n\alpha), \qquad 1\leq n\leq q-1.
$$
The hopping terms in the three Dirichlet Hamiltonians are
identical, so this proves
\eqref{eq:rational-irrational-sandwich}.

Consider now the eigenvalue counting functions for irrational 
$$
\widehat{\mathcal N}_q(x,E)
 :=
 \mathrm{Tr}\mathbf 1_{(-\infty,E]}
 \bigl(\widehat D_q(x)\bigr)
$$
and rational frequency
$$
\mathcal N_q^D(E)
 :=
 \mathrm{Tr}\mathbf 1_{(-\infty,E]}
 \bigl(D_q(0)\bigr).
$$
By the min-max principle,
\begin{equation}
\label{eq:counting-sandwich}
\widehat{\mathcal N}_q(x_q^-,E)
\geq
\mathcal N_q^D(E)
\geq
\widehat{\mathcal N}_q(x_q^+,E).
\end{equation}

By Corollary~3.6 of \cite{
jk2} we can estimate the difference between empirical DOS and the infinite-volume DOS measure
$$
\left|
 \frac1{q-1}\widehat{\mathcal N}_q(x,E)-N(E)
\right|
\leq
\frac{C}{q}
$$
with a universal constant that does not depend on $q$, $x$ or $E$. Combining this estimate with
\eqref{eq:counting-sandwich} proves
\eqref{eq:rational-Dirichlet-counting-discrepancy}.

\end{proof}

This comparison is particularly simple here, since we never come very close to the jump discontinuity;  we would have encountered the jump exactly at the site $0$ which is not part of the interval. It is valid also more generally as monotonicity up to finite-rank perturbation is enough at the cost of increasing the constant $C$.

A known consequence of this bound for the density of states combined with the Koksma inequality is the following large-deviations principle:
\begin{corollary}
For every bounded-variation function $g\in \mathrm{BV}(\RM)$,
\begin{equation}
\label{eq:rational-Dirichlet-Koksma}
\left|
\frac{1}{q-1} \sum_{i=1}^{q-1} g(\mu_i)
 -
 \int_{\RM}g(E')\,dN(E')
\right|
\leq
\frac{C}{q}\norm{g}_{\mathrm{BV}(\RM)}
\end{equation}    
with some universal constant $C$ and $\|g\|_{\mathrm{BV}(\mathbb R)}$ the total variation of $g$.
\end{corollary}
\begin{proof}
With the previous Lemma in place this can be proven exactly as \cite[Corollary~3.7]{
jk2}.
\end{proof}

\begin{lemma}
\label{lem:Dirichlet-LDT}
Let $\mu_1,...,\mu_{q-1}$ be the eigenvalues of $D_q$. For every $\epsilon>0$ and compact $K\subset \RM$ there exists some $Q_{K,\epsilon}>0$ such that
\begin{equation}
\label{eq:Dirichlet-punctured-determinant}
\frac{1}{q-1}
\log \prod_{i:\, i\neq k} \left|\mu_{i}-\mu_{k}\right|
\geq 
L(\mu_{k}) - \epsilon
\end{equation}
for every $\mu_k \in K$ and all $q\geq Q_{K,\epsilon}$.
\end{lemma}

\begin{proof}
As in \cite{jk2} we compare to the Lyapunov exponent via the Thouless formula 
$$L(E) := \int_{\mathbb{R}} \log \lvert E - E' \rvert \, dN(E')$$
by approximating the logarithm with a bounded-variation function 
$$
g_B(E')= \log \lvert E - E' \rvert_{[B_-,B_+]},\qquad 
|t|_{[B_{-},B_+]}
 :=
 \begin{cases}
 1, & |t|<B_{-},\\
 |t|, & B_{-}\leq |t|\leq B_+,\\
 1, & |t|>B_+
 \end{cases}
$$
for two cutoffs $B_- < 1  < B_+$. For the lower cutoff we will choose $B_-= \frac{\gamma}{2q}$ and for the upper cutoff $B_+=q$. With this choice we have
\begin{align*}
\frac{1}{q-1}\log\prod_{i:\, i\neq k} \left|\mu_{i}-\mu_{k}\right|&= \frac{1}{q-1}\sum_{i:\, i\neq k}\log\left|\mu_{i}-\mu_{k}\right|\\
&= \frac{1}{q-1}\sum_{i=1}^{q-1} g_B(\mu_{i})+  \frac{1}{q-1}\sum_{\substack{i:\, i\neq k\\ |\mu_{i}-\mu_{k}|>B_+}}\log\left|\mu_{i}-\mu_{k}\right|\\
&\geq \frac{1}{q-1}\sum_{i=1}^{q-1} g_B(\mu_{i}) \geq   \int_{\mathbb{R}} g_B(E')\, dN(E') - \frac{C}{q}\|g_B\|_{\mathrm{BV}(\RM)}
\end{align*}
where it is understood that we set $E=\mu_k$ in $g_B$. Here, we used that the product omits the eigenvalue $\mu_{k}$ and the next-closest eigenvalue is spaced at least $\frac{\gamma}{q}$ away from it, hence we can substitute $g_B$ below the higher cutoff $B_+$. Then we used that the contributions from $|\mu_i-\mu_k|>B_+$ are positive, so dropping them only makes the sum smaller. Finally we applied the Koksma inequality \eqref{eq:rational-Dirichlet-Koksma} to the bounded-variation function $g_B$ which satisfies
$$\|g_B\|_{\mathrm{BV}(\RM)}\leq 4 \log \frac{B_+}{B_-}.$$
It remains to compare the truncated integral to the Lyapunov exponent. We have
$$\int_{\mathbb{R}} g_B(E')\, dN(E') \geq L(E)-\int_{|E-E'|>B_+} \log|E-E'|\, dN(E')$$
since omitting the lower cutoff $B_-$ can only make the integral more negative. Put $R= \sup_{E\in K}|E|$ such that  
$$|E-E'|\leq |E|+|E'| \leq R+ |E'|.$$

We can then bound
$$\int_{|E-E'|>B_+} \log|E-E'|\, dN(E')\leq \int_{|E'|>B_+-R} \log(1+R+|E'|) dN(E').$$
As $|N^{-1}(t)-f(t)|\leq 2$, the logarithmic integrability \eqref{eq:log-integrability} of $f$ gives
\begin{align*}
\int_{\RM} \log(1+R+|E'|) dN(E') &\leq \int_{0}^1 \log(3+R+|f(t)|) dt\\
 &\leq \log(3+R) + \int_{0}^1 \log(1+|f(t)|) dt < \infty.
\end{align*}
Since $B_+=q\to \infty$ the integrability implies $$\int_{\mathbb{R}} g_B(E')\, dN(E') \geq L(E)-\frac{\epsilon}{2}$$
for all large enough $q$ uniformly for all $E\in K$, which together with $\frac{C}{q}\|g_B\|_{\mathrm{BV}(\RM)} \to 0$ completes the proof.
\end{proof}

This Lemma combined with the identity \eqref{eq:endpoint-cofactor-identity} finishes the proof of Lemma~\ref{lemma:endpoint_lyapunov_decay}.

\section{Quantitative spectral continuity}
\label{sec:continuity}
In this section we prove Theorem~\ref{th:quant_cont} and Theorem~\ref{th:lyapunov_controlled_approx}.

Let us briefly recall the usual residual argument that approximate eigenvectors give bounds for the distance to the spectrum as they imply largeness of the resolvent at the given energy. For possibly unbounded potentials the argument can be made quantitative as follows:
\begin{lemma}
\label{lemma:trial_vector}
Let $V,\tilde{V}$ be real multiplication operators on $\ell^2(\ZM)$ and let $\varphi$ be a normalized vector in $\mathrm{Dom}(V)$ with $\|(\Delta+V-E)\varphi\|\leq \epsilon$. If the resolvent-difference satisfies
$$\sup_{n\in \mathrm{supp}(\varphi)}\left|(V(n)+\imath)^{-1} - (\tilde{V}(n)+\imath)^{-1}\right| \leq \delta$$
with some $\delta>0$ so small that 
$$\delta (3+|E|+\epsilon) < \frac{1}{2}$$
then
\begin{equation}
\label{eq:spec_distance}\mathrm{dist}(E, \sigma(\Delta+\tilde{V}))\leq 2\epsilon +  2(3+ |E|+\epsilon)^2 \delta.
\end{equation}
\end{lemma}
\begin{proof}
Define the auxiliary vector $\tilde{\varphi}=(V+\imath)(\tilde{V}+\imath)^{-1}\varphi$ and note that by definition $\tilde{\varphi}\in \mathrm{Dom}(\tilde{V})$. Since $\|\Delta\|=2$ one has $$\|(V+\imath)\varphi\|=\|(\Delta+V-E)\varphi + (\imath+E-\Delta)\varphi\|\leq 3+ |E|+\epsilon$$ and then
$$\|\varphi-\tilde{\varphi}\|=\|\left((V+\imath)^{-1}-(\tilde{V}+\imath)^{-1}\right)(V+\imath)\varphi\| \leq \delta (3+|E|+\epsilon).$$
As by assumption $\delta (3+|E|+\epsilon)<\frac12$ we conclude $\|\tilde{\varphi}\|\geq\frac{1}{2}$ by the inverse triangle inequality.

Consider now the trivial identities
$$(V+\imath)\varphi = (\Delta+V-E)\varphi - (\Delta-E-\imath)\varphi$$
and
$$(V+\imath)\varphi = (\tilde{V}+\imath)\tilde{\varphi} = (\Delta+\tilde{V}-E)\tilde{\varphi}  - (\Delta-E-\imath)\tilde{\varphi}.$$
Subtracting them from one another gives
$$(\Delta+\tilde{V}-E)\tilde{\varphi}=(\Delta+V-E)\varphi - (\Delta-E-\imath)(\varphi-\tilde{\varphi}).$$
Therefore
$$\|(\Delta+ \tilde{V}-E)\tilde{\varphi}\|\leq \epsilon + (3+ |E|+\epsilon)^2 \delta$$
and since $\|\tilde{\varphi}\|\geq \frac12$ the normalization of $\tilde{\varphi}$ gives a trial eigenvector with
$$\|(\Delta+ \tilde{V}-E)\frac{\tilde{\varphi}}{\|\tilde{\varphi}\|}\|\leq 2\epsilon + 2(3+ |E|+\epsilon)^2 \delta.$$
This finishes the proof of \eqref{eq:spec_distance} as $\mathrm{dist}(E, \sigma(\Delta+\tilde{V})) = \|(\Delta+\tilde{V}-E)^{-1}\|^{-1}$.

\end{proof}

Based on this strategy we will first prove a quantitative criterion for the distance to the spectrum given a trial eigenvector. For any continued-fraction denominator $q$ we denote by $q^+$ the next-largest denominator.
\begin{proposition}
\label{prop:spectrum_error_quant}
Let $f$ be $\rho$-resolvent-H\"older continuous and let $K\subset \RM$ be a compact interval. There exist constants $C_{\rho,K}$ depending on $K$, $\rho$ and the Hölder constant such that the following holds: If there exists a normalized vector $\varphi\in \ell^2(\ZM)$ supported in an interval $[j-\ell,j+\ell]$ with $\ell < \frac{1}{10}q^+$ and
$$\|(H(x)-E)\varphi\|\leq \epsilon < 1$$
for some $x \notin \ZM + \alpha \ZM$ and $E\in K$ then
$$\mathrm{dist}(E, \Sigma_q) \leq 2\epsilon + C_{\rho,K} |\ell \Delta_q|^\rho$$
for sufficiently large $q$ (also depending only on $\rho$ and $K$ and the Hölder constant).
\end{proposition}

Note that this statement would be evident from Lemma~\ref{lemma:trial_vector} if $f$ itself was $\rho$-Hölder continuous as a function on $\TM$ as the scale of the second term is the expected deviation in operator norm between the irrational potential and its rational approximation on an interval of length $O(\ell)$. The main problem is therefore to deal with the jump discontinuity. 

The trial eigenvectors $\varphi$ will be constructed by truncating a given approximate eigenvector to an interval $[j-\ell,j+\ell]$ with the center $j$ optimized to minimize the truncation error. We therefore cannot exclude the possibility that on the support of $\varphi$ some phases $x+n \alpha$, $n\in [j-\ell,j+\ell]$ are very close to the jump discontinuity with $\|x+n\alpha\|\ll \|q\alpha\|$. Then direct comparison with $q$-periodic approximation must fail: Going $q$ steps to the left or right, the two potential values $f(x+n\alpha \pm q\alpha)$ sample $f$ on two different sides of the jump discontinuity, whereas a $q$-periodic approximation will necessarily sample $f$ on the same side both times, giving a mismatch of $O(f(0-)-f(0+))$ in the limit. 

To circumvent this problem we can use that on scales $\ell \ll q^+$ such points where $q$-periodicity must break down are rare; in fact if such a point exists on the support of $\varphi$ then it is unique. In the worst case, we can therefore use $\varphi$ not necessarily as a trial eigenvector for a $q$-periodic approximation, but at least it will still be a good trial eigenvector for a domain-wall model, which uses one $q$-periodic potential to the right of some jump at location $n$ and another $q$-periodic potential to its left. 

Let us now formalize this. If we want to minimize the phase mismatch between rational and irrational approximation at $j$ we use the irrational and rational approximations to sample respectively the points
$$x_n = \{x+n\alpha\}, \qquad y_n=\{x+j\alpha + (n-j)\frac{p}{q}\}.$$ We will say that $j$ is resonant if
\begin{equation}
\label{eq:minimum_dist}
\min_{n\in [j-\ell,j+\ell]} \|y_n\|= \min_{n\in [-\ell,\ell]} \|x+j\alpha + n\frac{p}{q}\| < \ell \Delta_q,
\end{equation}
i.e. if $y_n$ comes very close to the jump discontinuity. Otherwise we say $j$ is non-resonant. Throughout the comparison we assume $\ell< \frac{1}{10}q^+$, which implies $\ell\Delta_q< \frac{1}{10q}$.

In the non-resonant case we can compare directly to a rational approximation: 
\begin{lemma}
\label{lem:res_comp_nonresonant}
If $j$ is non-resonant and $\ell<\frac{1}{10}q^+$ then 
\begin{equation}
\label{eq:res_holder_comp1}
\sup_{n\in [j-\ell,j+\ell]}\left\| \frac{1}{f(x_n)+\imath} -  \frac{1}{f(y_n)+\imath}\right\|\leq C_\rho |\ell \Delta_q|^\rho
\end{equation}
for $C_\rho$ the Hölder constant of $(f+\imath)^{-1}$.
\end{lemma}
\begin{proof}
We say that two phases $x,y\in \TM$ are on the same side if the geodesic $[x,y]$ does not pass through $0$, otherwise we say $x,y$ are on opposite sides. If $x_n$ and $y_n$ are on the same side one has by the resolvent-Hölder continuity
$$\left\| \frac{1}{f(x_n)+\imath} -  \frac{1}{f(y_n)+\imath}\right\|\leq C_\rho |\ell \Delta_q|^\rho.$$
Note that $x_n$ and $y_n$ are on the same side for every $n\in [j-\ell,j+\ell]$ as $\|y_n\|\geq \ell \Delta_q$ and $\|x_n-y_n\|\leq \ell \Delta_q$.
\end{proof}
We consider now the case where $j$ is resonant. By $q$-periodicity the minimum in \eqref{eq:minimum_dist} is always attained by some $n_*\in \ZM$ such that $|n_*-j|<q$. If $j$ is resonant then we slightly modify $y_n$ and instead sample on an orbit which passes exactly through the jump discontinuity 
$$\tilde{y}_n=\{ y_n - y_{n_*}\}.$$
Consider then a regularized $q$-periodic approximation 
$$\tilde{V}^{q}_{j,c}(n) := \begin{cases}
c & \text{if } n-n_* \in q\ZM\\
f(\tilde{y}_n) & \text{otherwise}
\end{cases}$$
which replaces the values at the jump discontinuity by some real number $c\in \RM$. The idea is that $c$ will have to be sometimes close to $f(0+)$ and sometimes close to $f(1-)$ to get a good approximation. We therefore define for some $n_0\in \ZM$ domain-wall models of the form 
$$\tilde{V}^{(q)}_{j, n_0, c|d}(n)=\begin{cases}
\tilde{V}^{q}_{j,c}(n) & \text{if } n < n_0\\
\tilde{V}^{q}_{j,d}(n) & \text{if } n \geq n_0
\end{cases}.$$
\begin{lemma}
\label{lem:res_comp_resonant}
For any resonant $j\in \ZM$ and $\ell\in \NM$ with $\ell< \frac{1}{10}q^+$ there exist choices $n_0, c,d$ such that
\begin{equation}
\label{eq:res_holder_comp2}
\sup_{n\in [j-\ell,j+\ell]}\left\| \frac{1}{\tilde{V}^{(q)}_{j,n_0, c|d}(n)+\imath} -  \frac{1}{f(x_n)+\imath}\right\|\leq C_\rho |2\ell \Delta_q|^\rho.
\end{equation}
One can choose $c,d\in \{f(\ell \Delta_q), f(1- \ell \Delta_q)\}$.
\end{lemma}
\begin{proof}
The potentials are sampled at the points $x_n=x+n\alpha$ and $\tilde{y}_n$ with $$\|x_j-\tilde{y}_j\|= \min_{n\in [j-\ell,j+\ell]} \|y_n\|<\ell \Delta_q$$
and the phase mismatch therefore satisfies
$$\|x_n-\tilde{y}_n\|\leq 2\ell \Delta_q$$
for all $n\in [j-\ell,j+\ell]$. 

Note that $x_n$ and $\tilde{y}_n$ are on the same side for any $n\notin n_* + q\ZM \cap [j-\ell,j+\ell]$ as $\|\tilde{y}_n\|\geq \frac{1}{q}$ for all such $n$ and $\|x_n-\tilde{y}_n\|\leq 2\ell \Delta_q\leq \frac{1}{5q}$ due to $\Delta_q\leq \frac{1}{qq^+}$. Assume therefore the remaining case that $n$ is in the set
$$M=[j-\ell,j+\ell] \cap (n_*+q\ZM).$$
As $\tilde{y}_n=0$ for $n\in M$ we have $x_n\in [1-2\ell\Delta_q, 2\ell\Delta_q)$. There are then four possibilities:
\begin{enumerate}
    \item $x_{n} \in [0,2\ell \Delta_q)$ for each $n\in M$.
    \item $x_{n} \in [1-2\ell \Delta_q,1)$ for each $n\in M$.
    \item There exists some $n_0\in M$ such that $x_{n} \in (0,2\ell \Delta_q)$ for $n<n_0$ and $x_n \in (1-2\ell \Delta_q,1)$ for $n\geq n_0$.
    \item There exists some $n_0\in M$ such that $x_{n} \in (1-2\ell \Delta_q,1)$ for $n<n_0$ and $x_n \in (0,2\ell \Delta_q)$ for $n\geq n_0$.
\end{enumerate}
In the cases (3) and (4) such $n_0$ is unique as the phase $m\mapsto x_{n_0+mq}$ rotates by steps of $q\alpha$ and stays in $[1-2\ell\Delta_q, 2\ell\Delta_q]$, therefore cannot cross $0$ more than once.

Define now the two real numbers $$a=f(\ell \Delta_q), \quad b= f(1-\ell \Delta_q)$$ which will replace the possibly infinite $f(0+)$ and $f(1-)$ by finite values that are close enough for the required scale.  

For $n\notin M$ one has as before 
$$\left\| \frac{1}{f(x_n)+\imath} -  \frac{1}{f(\tilde{y}_n)+\imath}\right\|\leq C_\rho |2\ell \Delta_q|^\rho,$$
with $C_\rho$ the Hölder constant. 

For $n\in M$, if $x_{n} \in (0,2\ell\Delta_q)$ then 
$$
\left\| \frac{1}{a+\imath} -  \frac{1}{f(x_n)+\imath}\right\| \leq  C_\rho | 2\ell \Delta_q|^\rho$$
and if $x_{n}\in (1-2\ell\Delta_q, 1)$ one has
$$\left\| \frac{1}{b+\imath} -  \frac{1}{f(x_n)+\imath}\right\| \leq   C_\rho | 2\ell \Delta_q|^\rho.$$
If $j$ is resonant and $n_0$ either the unique crossover point or otherwise choosing arbitrarily $n_0\notin [j-\ell,j+\ell]$ and for $(c,d)$ the appropriate choice out of $\{(a,a),(b,b),(a,b), (b,a)\}$ one concludes \eqref{eq:res_holder_comp2}.

\end{proof}
This lemma together with the trial vector argument shows that for given $E\in \Sigma$ a suitable domain-wall Hamiltonian
$$H_{j,n_0,c|d}^{(q)}=\Delta+ \tilde{V}^{(q)}_{j,n_0,c|d}$$
must have some spectrum close to $E$. However, this is a priori unrelated to the periodic union spectrum $\Sigma_q$. We note first the technical result that interface eigenvalues are absent in our case, hence the spectra decouple:
\begin{lemma}
\label{lem:domain_wall_spec}
For any $c, d$ in $\RM$ and $n_0\in \ZM$ one has
$$\sigma(H^{(q)}_{j, n_0, c|d})= \sigma(H^{(q)}_{j,c})\cup \sigma(H^{(q)}_{j,d})$$
with 
$$H_{j,c}^{(q)}=\Delta+ \tilde{V}^{(q)}_{j,c}.$$
\end{lemma}
The proof will be given at the end of the section since it is a more general fact about domain walls of $q$-periodic models. The spectra of the $q$-periodic models themselves are still not necessarily part of the union spectrum $\Sigma_q$ as $H^{(q)}_{j,c}$ is only a rational approximation if $c=f(0+)$. But they are close:
\begin{lemma}
\label{lem:res_spec_comp}
Fix some compact interval $K\subset \RM$. If $\ell < \frac{1}{10}q^+$, $j\in \ZM$ is resonant, $c\in \{f(\ell\Delta_q), f(1-\ell\Delta_q)\}$ and $E\in \sigma(H_{j,c}^{(q)})\cap K$ then
$$\mathrm{dist}(E,\Sigma_q) \leq C_{K,\rho}|\ell\Delta_q|^\rho$$ for every large enough $q$ and with a constant depending only on $K$ and the $\rho$-Hölder constant of $f$. \end{lemma}
\begin{proof}
Recall that $H_{j,c}^{(q)}$ samples at the $q$-periodic orbit $\tilde{y}_n$ with $n\in \ZM$ and which passes through the phase $0$. We compare to the $q$-periodic orbit $\tilde{y}_n\pm \ell\Delta_q$ with positive sign if $c=f(\ell\Delta_q)$ and negative if $c=f(1-\ell\Delta_q)$. Denote the corresponding $q$-periodic potential sampling the values $f(\tilde{y}_n\pm \ell\Delta_q)$ by $W_\pm$. The spectrum $\sigma(\Delta+W_\pm)\subset \Sigma_q$ of this rational approximation is by definition included in the periodic union spectrum. Then resolvent-Hölder continuity gives
$$\|(\tilde{V}^{(q)}_{j,c}(n)+\imath)^{-1}- (W_\pm(n)+\imath)^{-1}\| \leq C_\rho |\ell\Delta_q|^\rho, \qquad \forall n\in \ZM$$
as each value of the potential compares either $f(\tilde{y}_n\pm \ell\Delta_q)$ to $f(\tilde{y}_n)$ and both are on the same side, or it compares $c=f(\ell\Delta_q)$ to $f(0+)$ or it compares $c=f(1-\ell\Delta_q)$ to $f(1-)$.

For a normalized trial eigenvector $\psi$ of $H^{(q)}_{j,c}$ with error $\epsilon>0$ Lemma~\ref{lemma:trial_vector} implies 
$$\mathrm{dist}(E, \Sigma_q)\leq 2\epsilon+2(3+|E|+\epsilon)^2 C_\rho (\ell\Delta_q)^\rho$$
whenever $(3+|E|+\epsilon)C_\rho(\ell\Delta_q)^\rho< \frac{1}{2}$. Letting $\epsilon$ go to $0$ gives the conclusion.
\end{proof}

\begin{proof}(of Proposition~\ref{prop:spectrum_error_quant})
Let $\varphi$ be a trial eigenvector constructed for $H(x)$. If $j$ is non-resonant we compare with the rational approximation on the orbit $y_n$, otherwise we compare with the suitable $H^{(q)}_{j,n_0,c|d}$ on an orbit $\tilde{y}_n$. In either case let $\tilde{V}^{(q)}$ be the corresponding $q$-periodic or domain-wall potential. 
Assuming $\ell< \frac{1}{10}q^+$ the residual of the trial eigenvector is in either case using Lemma~\ref{lem:res_comp_nonresonant} respectively Lemma~\ref{lem:res_comp_resonant} determined by
\begin{align*}
\sup_{n\in [j-\ell,j+\ell]}\left|(\tilde{V}^{(q)}(n)+\imath)^{-1}-(f(x_n)+\imath)^{-1}\right| &\leq C_\rho |2\ell \Delta_q|^\rho
\end{align*} 
such that Lemma~\ref{lemma:trial_vector} implies
$$\mathrm{dist}(E, \sigma(\Delta+\tilde{V}^{(q)}))\leq 2 \epsilon + C_{K,\rho} |\ell \Delta_q|^\rho.$$
In the non-resonant case this also bounds $\mathrm{dist}(E,\Sigma_q)$ immediately. In the resonant case Lemma~\ref{lem:domain_wall_spec} and Lemma~\ref{lem:res_spec_comp} give $$\mathrm{dist}(E, \Sigma_q)\leq \mathrm{dist}(E, \sigma(\Delta+\tilde{V}^{(q)}))+ C_{K,\rho} |\ell\Delta_q|^{\rho}$$
hence the additional error is again of the same size.
\end{proof}

Finally, we supply an elementary proof for Lemma~\ref{lem:domain_wall_spec}, which is a special case of the following:

\begin{lemma}
\label{lem:one_site_domain_wall}
Fix some function $v:\{1,...,q-1\}\to \RM$ and define the $q$-periodic potential $$V_a:\ZM\to \RM, \quad V_a(n)=\begin{cases}
v_r & n=mq+r, \quad r\neq 0\\
a & n=mq.
\end{cases}
$$
For $a,b\in \RM$ and the domain-wall potential
$$V_{a|b}(n)=\begin{cases}
v_r, \qquad &\text{if }  n=mq+r, \quad r\neq 0\\
a, \qquad &\text{if }    n=mq, \quad m < 0\\
b, \qquad &\text{if }    n=mq, \quad m\geq 0
\end{cases}$$
one has
$$\sigma(\Delta+V_{a|b})=\sigma(\Delta+V_a)\cup \sigma(\Delta+V_b).$$
\end{lemma}
\begin{proof}
The inclusion $$\sigma(\Delta+V_{a|b})\supset \sigma(\Delta+V_a)\cup \sigma(\Delta+V_b)$$ is standard using Weyl sequences. By another standard argument one has
 $$\sigma_{\mathrm{ess}}(\Delta+V_{a|b})= \sigma_{\mathrm{ess}}(\Delta+V_a)\cup \sigma_{\mathrm{ess}}(\Delta+V_b).$$
Suppose for a contradiction that $E\in \sigma(\Delta+V_{a|b})$ is not in $\sigma(\Delta+V_a)\cup \sigma(\Delta+V_b)$, in particular it must then be an isolated eigenvalue. Let $\psi$ be the associated eigenvector and consider the restrictions $\psi_m\in \CM^{q-1}$ of $\psi$ to the cells $\{mq+1,...,(m+1)q-1\}$ and $\xi_m=\psi(mq)$. Let $D$ be the Schrödinger operator with Dirichlet boundary conditions on the common cell $\{1,...,q-1\}$ and assume first $E\notin \sigma(D)$. As the eigenvalue equation gives 
\begin{equation}
\label{eq:block_ev_equation}
(E-D)\psi_m =\xi_m\delta_1 + \xi_{m+1} \delta_{q-1}
\end{equation}
we can multiply with the resolvent $R_E=(E-D)^{-1}$ to obtain 
$$\psi_m = \xi_m (R_E\delta_1) + \xi_{m+1} (R_E \delta_{q-1}).$$ Indeed, this is just the Poisson formula written in vector form. The eigenvalue equation at $mq$ written in a similar form reads
$$\langle \delta_1, \psi_m\rangle + \langle \delta_{q-1}, \psi_{m-1}\rangle + c_m \xi_m = E \xi_m$$
with $c_m\in \{a,b\}$, and substituting the expressions for $\psi_m$ we can rewrite it as
$$ \xi_m \langle \delta_1, R_E\delta_1\rangle + \xi_{m+1} \langle \delta_1,R_E \delta_{q-1}\rangle +\xi_{m-1} \langle  \delta_{q-1}, R_E\delta_1\rangle + \xi_{m} \langle \delta_{q-1},R_E \delta_{q-1}\rangle = (E-c_m) \xi_m$$
which is the difference equation
$$\xi_{m-1} + \xi_{m+1} = F_m(E) \xi_m, \qquad F_m(E)=\frac{E-c_m- \langle \delta_1, R_E\delta_1\rangle -\langle \delta_{q-1}, R_E\delta_{q-1}\rangle}{\langle \delta_{q-1}, R_E \delta_1\rangle}.$$
We used that by Cramer's rule $|\langle \delta_{q-1}, R_E \delta_1\rangle|=|\det(E-D)^{-1}|$, which is non-zero.

A normalizable eigenfunction can therefore only have the functional form
\[
\xi_m=A\lambda_a^m\quad(m\leq-1),
\qquad
\xi_m=B\lambda_b^{-m}\quad(m\geq1)
\]
with some $|\lambda_a|,|\lambda_b|>1$ that are fixed uniquely by the characteristic equations
$$\lambda_a + \lambda_a^{-1}=F_{-1}(E), \qquad \lambda_b + \lambda_b^{-1}=F_1(E).$$
The equations at \(m=-1\) and \(m=1\) together force
\[
\xi_0=A=B
\]
and the equation at \(m=0\) then becomes
\[
A\lambda_a^{-1}+A\lambda_b^{-1}
=
F_0(E)A=F_1(E)A
=
A\lambda_b+A\lambda_b^{-1}.
\]
For $A\neq 0$ this gives $\lambda^{-1}_a=\lambda_b$ which immediately contradicts $|\lambda_a|>1$ and $|\lambda_b|>1$. If conversely $A=0$, then $\psi$ vanishes on every site in $q\ZM$, hence $E\in \sigma(D)$, which is a contradiction to the assumption made above.

It remains to consider the case \(E\in\sigma(D)\). Let $u$ be a real eigenvector of $D$ with eigenvalue $E$ and consider the scalar product of \eqref{eq:block_ev_equation} with $u$ to get
$$0 = u(1) \xi_m + u(q-1)\xi_{m+1}.$$ As both $u(1)$ and $u(q-1)$ are non-vanishing, this implies $\xi_m = c \lambda^m$ is a geometric sequence with $\lambda\neq 0$. Thus, $\psi$ can only be an eigenfunction if $c=0$, i.e. if $\psi$ vanishes on $q\ZM$. However, that would imply that the $q$-periodic extension of $D$ to $\ell^2(\ZM)$ has $\psi$ as an $\ell^2$-normalizable eigenfunction, which is not possible (a periodic Schrödinger operator cannot have point spectrum).
\end{proof}

\begin{remark}
In this section we assumed that $f$ has one (periodically repeated) jump discontinuity. A similar analysis can be performed for piecewise continuous $f$, provided that the discontinuities are well-separated, e.g. satisfy a relative Diophantine condition. Without such a separation condition multiple discontinuities could be resonant on the support of a given trial eigenvector, such that one would need more complicated local models than the above domain walls and absence of interface spectrum should not be expected to generalize to that situation. Hence it is to be expected that Proposition~\ref{prop:spectrum_error_quant} and ultimately also Theorem~\ref{th:quant_cont} and Theorem~\ref{th:lyapunov_controlled_approx} do not immediately generalize to the case of sampling functions which have more than one jump discontinuity. $\diamond$
\end{remark}

\subsection{Proof of \texorpdfstring{Theorem~\ref{th:quant_cont}}{Theorem 1.3}}
For a generic approximate eigenvector without any decay assumptions there exists an optimal choice of truncation \cite{avms}: 

\begin{lemma}
\label{lem:trial_const}
For $H=\Delta+V$ with $V$ a real multiplication operator, any $\epsilon>0$ and $\ell\in \NM$, $E\in \sigma(H)$ there exists some $j\in \ZM$ and a normalized trial eigenvector $\varphi\in \ell^2(\ZM)$ supported in $[j-\ell,j+\ell]$  such that
$$\|(H-E)\varphi\|\leq C (\epsilon+ \ell^{-1})$$
with some universal constant $C$.
\end{lemma}
\begin{proof}
Let $\epsilon>0$ then one can always pick a normalized trial eigenvector $\phi\in \ell^2(\ZM)\cap\mathrm{Dom}(H)$ with residual 
$$\|(H-E)\phi\|< \epsilon.$$
We truncate it by multiplying it with the soft cutoff $$g_{j,\ell}(n)=\left(1- \frac{|n-j|}{\ell}\right)_+.$$ 
Since $|g_{j,\ell}(n+1)-g_{j,\ell}(n)|\leq \ell^{-1}$ one can also see that 
$$\sum_{j}\|[\Delta, g_{j,\ell}]\phi\|^2 \leq 4\ell^{-2} \sum_{j}\|\chi_{j,\ell}\phi\|^2 = 4(2\ell+3)\ell^{-2}\|\phi\|^2\leq C \ell^{-1}$$
with $\chi_{j,\ell}$ the characteristic function of the interval $[j-\ell-1,j+\ell+1]$ where the middle equality used that each $|\phi(n)|^2$ for fixed $n$ appears in the sum $(2\ell+3)$ times.

We can use this to commute the approximate eigenfunction past $H$ and bound
\begin{align*}
\left(\sum_{j}\|(H-E)g_{j,\ell}\phi\|^2\right)^{\frac12} &\leq\left(\sum_{j}\|g_{j,\ell}(H-E)\phi\|^2\right)^{\frac12}+\left(\sum_{j}\|[\Delta,g_{j,\ell}]\phi\|^2\right)^{\frac12}\\
&\leq  (S_\ell^\frac12 \epsilon + \sqrt{C}\ell^{-1/2})
\end{align*}
with
$$S_\ell = \sum_{j}\|g_{j,\ell}\phi\|^2.$$
Noting that $S_\ell\asymp \ell$ one can also put the inequality in the form
$$\sum_{j}\|(H-E)g_{j,\ell}\phi\|^2 \leq  C' \left(\epsilon + \ell^{-1}\right)^2 S_\ell$$
with a larger but still universal constant $C'$.
Terms with $g_{j,\ell}\phi = 0$ vanish on both sides. If every $j$ with $g_{j,\ell}\phi \neq 0$ would satisfy $$\|(H-E)g_{j,\ell}\phi\|^2 >  C' \left(\epsilon + \ell^{-1}\right)^2 \|g_{j,\ell}\phi\|^2$$
then summing over $j$ would give a contradiction. Hence there is some $j$ with
$$\frac{\|(H-E)g_{j,\ell}\phi\|}{\|g_{j,\ell}\phi\|} \leq  \sqrt{C'} \left(\epsilon + \ell^{-1}\right).$$
The normalization of this $g_{j,\ell}\phi$  is therefore the desired trial vector $\varphi$.
\end{proof}

We now merely have to plug this choice into Proposition~\ref{prop:spectrum_error_quant}:

\begin{proof}(of Theorem~\ref{th:quant_cont})
Choose any generic phase $x \notin \ZM+\alpha\ZM$ such that $\sigma(H(x))=\Sigma$ and let $E\in \Sigma\cap K$.
For $\epsilon>0$ and $\ell>0$ to be determined let $\varphi$ be a trial eigenvector constructed as in Lemma~\ref{lem:trial_const} for $H=H(x)$ and the given energy $E$. 
Assuming $\ell< \frac{1}{10}q^+$ we can apply Proposition~\ref{prop:spectrum_error_quant} to get
$$\mathrm{dist}(E, \Sigma_q)\leq 2 C(\epsilon+\ell^{-1}) + C_{K,\rho}|\ell \Delta_q|^{\rho}.$$ 

For the optimal choice $\ell = c|\Delta_q|^{-\frac{\rho}{1+\rho}}+O(1)$, which makes the leading order errors the same size, one gets
$$\mathrm{dist}(E, \Sigma_q)\leq 2C\epsilon + C_{K,\rho}|\Delta_q|^{\frac{\rho}{1+\rho}}.$$ 
with $C_{K,\rho}$ depending only on $K$ and $C_\rho$.
This choice also satisfies the required  $\ell< \frac{1}{10}q^+$ for $c<\frac{1}{10\sqrt{2}}$, as $\Delta_q>\frac{1}{q(q+q^+)}$ implies $\ell< c(2(q^+)^2)^{\frac{\rho}{1+\rho}}\leq c \sqrt{2}q^+$. Taking the limit $\epsilon\to 0$ finishes the proof.
\end{proof}

\subsection{The case of positive Lyapunov exponent}

In this section we prove Theorem~\ref{th:lyapunov_controlled_approx}. We fix some Diophantine $\alpha$ and assume as in the statement of the theorem that $J$ is an interval on which the Lyapunov exponent $L(E)$ is continuous and has a strictly positive lower bound $L_J>0$.

Define the transfer matrix
$$A_E(x)=\begin{pmatrix} E-f(x) & -1 \\ 1 & 0\end{pmatrix}.$$
We assume throughout the conditions of Theorem~\ref{th:lyapunov_controlled_approx}. Since $f$ is piecewise continuous and satisfies the logarithmic integrability \eqref{eq:log-integrability} the Thouless formula gives the Lyapunov exponent of the transfer matrix.  To connect to the resolvent and make the transfer matrices bounded we renormalize to
$$B_E(x)=r(x) A_E(x) = \begin{pmatrix} (E+\imath)r(x)-1 & -r(x) \\ r(x) & 0\end{pmatrix}, \qquad r=(f+\imath)^{-1}.$$
The Lyapunov exponent for $B_{E,n}(x)=B_E(x+(n-1)\alpha)... B_E(x)$ is equal to
$$\lambda(E) = L(E) + \int_0^1 \log|r(x)| dx =: L(E) + L_f.$$

The first Lemma that we  need is uniform-in-phase control over the growth rate of both factors. 

Since we allow $-\infty$ limit at $0,$ we need to only involve $x$ whose orbit does not hit a discontinuity to insure $\sigma(H(x))=\Sigma$ where $\Sigma$ is the almost sure spectrum. Call a phase $x$ admissible if its $\alpha$-orbit $x+\ZM\alpha$ does not intersect the integers. The following is immediate from \cite[Theorem 3.1]{JitomirskayaMavi17} applied to the bounded subadditive cocycles $\log\|B_{E,n}\|$:
\begin{lemma}
\label{lemma:transfer_matrix_growth}
For every small enough $\epsilon>0$ there is some $n_0$ such that
$$\|B_{E,n}(x)\| \leq e^{n(\lambda(E)+\epsilon)}$$
for all $n\geq n_0$, $E\in J$ and all admissible phases $x$. Moreover, it can be chosen large enough that also
$$\prod_{j=0}^{n-1} |r(x+j\alpha)| \leq e^{n(L_f +\epsilon)}$$
for all $n\geq n_0$.
\end{lemma}
\begin{proof}
Let us first argue that $h_n(x)=\log\|B_{E,n}(x)\|$ is bounded for every fixed $n$. Since the norm $\|B_E(x)\|$ is bounded we only have to assert that $\inf_{x\in \TM}\|B_{E,n}(x)\|>0$. These functions are continuous away from the finitely many points $-j\alpha$, $0\leq j<n$. At each of those exceptional points both left and right limits exist and are non-zero matrices $B_{E,n}(-j\alpha\pm)$. To see this, note that the only non-invertible factor that can appear in a product $B_{E,n}(x\pm)=B_E(x+(n-1)\alpha\pm)...B_E(x\pm)$ is the singular matrix $$\begin{pmatrix}
-1 & 0 \\ 0 & 0
\end{pmatrix}$$ which can occur as $B_E(0+)$ or $B_E(1-)$ if $f$ tends to $+\infty$ or $-\infty$ as left or right limit. This factor can appear at most once since $0$ and $1$ are identified and irrational rotations are aperiodic, hence the product has at least rank $1$. By compactness and piecewise continuity, $x\mapsto \|B_{E,n}(x)\|$ has strictly positive lower bound for every $x$. 

We can therefore apply \cite[Theorem 3.1]{JitomirskayaMavi17} to see that one has for every $E\in J$ a neighborhood $B_\delta(E)$ such that $$\frac{1}{n}\log\|B_{E',n}(x)\|\leq \lambda(E)+ \epsilon \max(1, \|h_1\|_\infty)$$
for every $x\in \TM$, $E' \in B_\delta(E)$ and $n$ large enough. The continuity of $L(E)$ also implies continuity of $\lambda(E)$, hence absorbing the factor $\|h_1\|_\infty$ and using a finite covering argument implies that one can also pick some $n_0$ such that
$$\frac{1}{n}\log\|B_{E,n}(x)\|\leq \lambda(E)+ \epsilon$$
for every $E\in J$ and $x\in \TM$. 

The secondary claim follows from convergence
$$\frac{1}{n}\sum_{j=0}^{n-1}\max(\log|r(x+j\alpha)|,-M)\to \int_\TM \max(\log|r(y)|,-M) dy$$
which is uniform in $x$ for every cutoff $M>0$ since the truncated function is bounded and Riemann
integrable, and averages of such functions under an irrational rotation converge uniformly in the phase $x$. As $\int_\TM \log|r(y)| dy = L_f < \infty$ is integrable by assumption we can pick $M$ large enough such that the right-hand side is at most $L_f+\frac{\epsilon}{2}$, which then implies
$$\frac{1}{n}\sum_{j=0}^{n-1}\log|r(x+j\alpha)| \leq \frac{1}{n}\sum_{j=0}^{n-1}\max(\log|r(x+j\alpha)|,-M) \leq L_f +\epsilon$$
for large $n$ as the cutoff only makes the sum larger.
\end{proof}

The next step is to construct a set of good phases where one also has a lower bound on the transfer matrices.
\begin{lemma}
For every $\epsilon>0$ there exists some $n_0\in \NM$ such that for every $E\in J$ and $n\geq n_0$ there is an interval $I_{E,n}\subset \TM$ of length $|I_{E,n}|>e^{-\epsilon n/\rho}$ such that
$$\|B_{E,n}(x)\|\geq e^{n (\lambda(E)-\epsilon)}$$
for every admissible phase $x\in I_{E,n}$.
\end{lemma}
\begin{proof} Put $\delta_n=e^{-\varepsilon n/\rho}$. By subadditivity,
\[
\int_{\mathbb T}\log\|B_{E,n}(x)\|\,dx\ge n\lambda(E),
\]
so there is an admissible phase $x_0$ with
$\|B_{E,n}(x_0)\|\ge e^{n\lambda(E)}$.

Since $\alpha$ is Diophantine, the points
$\{-j\alpha:0\le j<n\}$ are separated by more than
$4\delta_n$ for all sufficiently large $n$.
Thus the continuity interval containing $x_0$ contains
an interval $I_{E,n}$ of length $2\delta_n$ with $x_0$
as one endpoint.

Lemma~4.10, with shorter products absorbed into a constant,
and the uniform $\rho$-H\"older continuity of $B_E$ on $(0,1)$ give,
by telescoping,
$$
\|B_{E,n}(x)-B_{E,n}(x_0)\|
\le C_\varepsilon n|x-x_0|^\rho   e^{n(\lambda(E)+\varepsilon/2)}
\le C'_\varepsilon n e^{n(\lambda(E)-\varepsilon/2)}
$$
for every admissible $x\in I_{E,n}$. Consequently,
$
\|B_{E,n}(x)\|
\ge e^{n\lambda(E)}
       \bigl(1-C'_\varepsilon n e^{-\varepsilon n/2}\bigr)
\ge e^{n(\lambda(E)-\varepsilon)}
$
for all sufficiently large $n$, uniformly in $E\in J$.
\end{proof}
The entries of the matrix $B_{E,n}$ are Dirichlet determinants up to factors of $r(x)$. If the norm of $\|B_{E,n}(x)\|$ is large then one of those determinants must also be comparably large. Via Cramer's rule this gives:
\begin{corollary}
\label{cor:decay_greensfunction}
Suppose that $0 < \epsilon< L_J/100$, then there exists some $n_0$ (uniformly for $E\in J$) such that for every $n\geq n_0$ and every admissible $x\in I_{E,n}$ there exists an integer interval $[n_1,n_2]\subset [0,n-1]$ with length $n-2\leq n_2-n_1+1 \leq n$ such that
$$|G^{E,x}_{[n_1,n_2]}(n_1,m)|\leq e^{-\frac{L_J}{3} n}$$
and
$$|G^{E,x}_{[n_1,n_2]}(m,n_2)|\leq e^{-\frac{L_J}{3} n}$$
for $m$ with $|m-\frac{n_1+n_2}{2}|\leq 2$, i.e. $m$ in the middle of the interval. 
\end{corollary}
\begin{proof}
Recall that the transfer matrix $A_{E,n}$ has as its entries the Dirichlet determinants
$$A_{E,n}(x)=(-1)^n	\begin{pmatrix}
\det(H_{[0,n-1]}(x)-E) & \det(H_{[1,n-1]}(x)-E)\\
-\det(H_{[0,n-2]}(x)-E) & -\det(H_{[1,n-2]}(x)-E)
\end{pmatrix}.
$$
In terms of the renormalized Dirichlet determinants
$$
D_{[n_1,n_2]}(x)
=\det(H_{[n_1,n_2]}(x)-E)\prod_{j=n_1}^{n_2}r(x+j\alpha),
$$
one therefore has 
$$
B_{E,n}(x)=(-1)^n
\begin{pmatrix}
D_{[0,n-1]}(x)
&
r(x)D_{[1,n-1]}(x)
\\
-r(x+(n-1)\alpha)D_{[0,n-2]}(x)
&
-r(x)r(x+(n-1)\alpha)D_{[1,n-2]}(x)
\end{pmatrix}.
$$
For $x\in I_{E,n}$ one has $\|B_{E,n}\|\geq e^{n(\lambda(E)-3\epsilon)}$ and therefore at least one entry of this matrix must have absolute value of at least $\frac{1}{4}e^{n(\lambda(E)-3\epsilon)}$. Let $[n_1,n_2]\subset [0,n-1]$ be the corresponding integer interval. Since $|r(x)|\leq 1$ one then also has $|D_{[n_1,n_2]}(x)|\geq \frac{1}{4}e^{n(\lambda(E)-3\epsilon)}$. For any $m\in [n_1,n_2]$ for which $n_2-m$ and $m-n_1$ are large enough, Cramer's rule for that interval gives
\begin{align*}
|G^E_{[n_1,n_2]}(n_1,m)| &= \frac{|D_{[m+1,n_2]}(x)|}{|D_{[n_1,n_2]}(x)|}
\prod_{j=n_1}^{m}|r(x+j\alpha)|\\
&\leq 4 e^{(n_2-m)(\lambda(E)+\epsilon)-n(\lambda(E)-3\epsilon)+(m-n_1+1)(L_f+\epsilon)}\\
&\leq 4 e^{-(m-n_1) L(E) + 2|\lambda(E)| + (1+4 n_2-4 n_1)\epsilon}.
\end{align*}
For $n$ large enough and $m$ close to the midpoint $|\frac{n_1+n_2}{2}-m|\leq 2$ we therefore have
$$|G^E_{[n_1,n_2]}(n_1,m)| \leq e^{-\frac{n}{2} L_J+5n\epsilon}\leq e^{-\frac{n}{3} L_J}$$
with the last inequality holding for small enough $\epsilon$ such as $\epsilon<L_J/100$. The same considerations apply to the right endpoint.
\end{proof}

We can use this to show that any generalized eigenfunction must be exponentially small at certain places and use that to construct short trial eigenvectors.

\begin{proof}(of Theorem~\ref{th:lyapunov_controlled_approx})
We assume that $\alpha$ is Diophantine with
$$\|m\alpha\|\geq c_\alpha |m|^{-\tau}$$
for all $m\neq 0$ and some $\tau>1$. For any admissible $x$ one has $\sigma(H(x))=\Sigma.$ By Schnol's Theorem (e.g. \cite{df}) and by redefining the phase and normalization, for spectrally almost every energy $H(x)$ has a generalized eigenvector 
$\varphi$ with
$$\varphi(0)=1, \qquad |\varphi(j)|\leq 1 + |j|.$$
For any given $E\in J \cap \Sigma$ we can pick now admissible phases $x_n$ and some such generalized eigenfunctions $\varphi_n$ with energies $E_n\in J$ and $|E-E_n|\leq e^{-an}$ for some constant $a>0$. Without loss of generality they are assumed to be shifted and normalized such that $\varphi_n(0)=1$ and $|\varphi_n(j)|\leq 1 + |j|$.

For any given $n$ let $Q_n$ be the first continued fraction denominator such that $Q_n > 100 e^{\epsilon n/\rho}$. By the Diophantine condition the size of $Q_n$ is at most
$$Q_n \leq  C_\alpha e^{\tau \epsilon n/\rho},$$
thus exponential in $n$. On a $Q_n$-step trajectory from any starting phase the largest gap is bounded by $\frac{3}{Q_n}$, hence within $Q_n$ steps in either direction from $j$ one is guaranteed to reach the interval $I_{E,n}$. We conclude that there are two good intervals of length between $n-2$ and $n$ as in Corollary~\ref{cor:decay_greensfunction} contained within $[-Q_n-n,-n]$ and $[1, Q_n]$ respectively. As a generalized eigenfunction $\varphi_n$ satisfies the Poisson formula
$$\varphi_n(m) = -G^{E_n,x_n}_{[n_1,n_2]}(n_1,m)\varphi(n_1-1)- G^{E_n,x_n}_{[n_1,n_2]}(m,n_2)\varphi(n_2+1).$$ Let $m_1, m_2$ be the midpoints of the intervals we just constructed. Truncating $\varphi_n$ at the midpoints of those intervals one has
\begin{align*}
\|(H(x_n)-E_n)\varphi_n \chi_{[m_1,m_2]}\|^2&= |\varphi_n(m_1-1)|^2+|\varphi_n(m_1)|^2+|\varphi_n(m_2)|^2+|\varphi_n(m_2+1)|^2 \\
&\leq e^{-2\frac{L_J}{3}n} C (1+Q_n)^2.
\end{align*}
Since $\|\varphi_n \chi_{[m_1,m_2]}\|>1$ the corresponding normalized trial vector has a properly smaller residual
\begin{equation}
\label{eq:error_residual}
\|(H(x_n)-E_n)\frac{\varphi_n \chi_{[m_1,m_2]}}{\|\varphi_n \chi_{[m_1,m_2]}\|}\| \leq C e^{-\frac{L_J}{3}n}(1+Q_n)
\end{equation}
for all large enough $n$. 

Since we have trial eigenvectors supported in $[-Q_n-n,Q_n+n]$ by Proposition~\ref{prop:spectrum_error_quant} we may conclude
\begin{equation}
\label{eq:error_spectrum}
\mathrm{dist}(E_n, \Sigma_q)\leq C e^{-\frac{L_J}{3}n} (1+Q_n) + C_{\rho,J} |(Q_n+n) \Delta_q|^\rho
\end{equation}
with a constant that depends only only on $\rho$ and $J$ (assuming $n$ and $Q_n$ are chosen in such a way that the Proposition can be applied).

For the given $0< s <\rho$ choose $\epsilon>0$ small enough such that $\frac{6s\tau\epsilon}{L_J \rho}< \min(\frac{s}{2}, 1-\frac{s}{\rho}, \frac{1}{2})$. For any continued fraction denominator $q$ we then pick $n=\lceil \frac{6s}{L_J} \log(\frac{1}{\Delta_q}) \rceil$ which makes all terms on the right-hand side of \eqref{eq:error_spectrum} smaller than 
$$\mathrm{dist}(E_n, \Sigma_q)\leq C_{\epsilon,s,J} |\Delta_q|^s$$
for sufficiently large $q$.

Finally, note that with this choice $\Delta_q \geq \frac{1}{2(q^+)^2}$ implies  $$\frac{Q_n}{q^+}\leq 2^{\frac{6s\epsilon\tau}{L_J \rho}} (q^+)^{\frac{12s\epsilon\tau}{L_J \rho}-1}\to 0$$ with $q\to \infty$, hence the support of $\varphi$ is smaller than $\frac{1}{10}q^+$ for large enough $q$, as was required to apply Proposition~\ref{prop:spectrum_error_quant}. Finally, choosing $a=\frac{1}{6}L_J$ makes $$|\mathrm{dist}(E_n,\Sigma_q)-\mathrm{dist}(E,\Sigma_q)|\leq e^{-an}\leq |\Delta_q|^s$$ another error of the same order.
\end{proof}


\begin{thebibliography}{99}


\bibitem{avilaarc1} A.~Avila,
{\it Almost reducibility and absolute continuity I},
preprint, arXiv:1006.0704 (2010).

\bibitem{avilaarc2} A.~Avila,
{\it KAM, Lyapunov exponents, and the spectral dichotomy for typical
one-frequency Schr\"odinger operators},
preprint, arXiv:2307.11071 (2023).

\bibitem{aj} A.~Avila and S.~Jitomirskaya,
{\it The Ten Martini Problem},
Annals of Mathematics (2) \textbf{170} (2009), no.~1, 303--342.

\bibitem{ak} A.~Avila and R.~Krikorian,
{\it Reducibility or nonuniform hyperbolicity for quasiperiodic
Schr\"odinger cocycles},
Annals of Mathematics (2) \textbf{164} (2006), no.~3, 911--940.

\bibitem{akmon}
A.~Avila and R.~Krikorian,
Monotonic cocycles,
\emph{Invent. Math.} \textbf{202} (2015), 271--331.



\bibitem{ayz} A.~Avila, J.~You, and Q.~Zhou,
{\it Dry Ten Martini Problem in the non-critical case},
preprint, arXiv:2306.16254 (2023).

\bibitem{ayz2} A.~Avila, J.~You, and Q.~Zhou, Preprint.

\bibitem{as} J.~E.~Avron and B.~Simon,
{\it Almost periodic Schr\"odinger operators II.
The integrated density of states},
Duke Mathematical Journal \textbf{50} (1983),
no.~1, 369--391.

\bibitem{avms} J.~E.~Avron, P.~H.~M.~van Mouche, and B.~Simon,
{\it On the measure of the spectrum for the almost Mathieu operator},
Communications in Mathematical Physics \textbf{132} (1990), 103-118.

\bibitem{bbl} R.~Band, S.~Beckus, and R.~Loewy,
{\it The Dry Ten Martini Problem for Sturmian Hamiltonians},
preprint, arXiv:2402.16703 (2024).

\bibitem{band-sofer} R.~Band and G.~Sofer,
{\it Johnson--Schwartzman gap labeling for metric and discrete
decorated graphs},
Journal of the London Mathematical Society \textbf{114} (2026),
no.~5, e70748.

\bibitem{bjk} S.~Becker, S.~Jitomirskaya, and I.~Krasovsky,
{\it Critical almost Mathieu operator: hidden singularity, gap continuity,
and the Hausdorff dimension of the spectrum},
preprint, \href{https://arxiv.org/abs/1909.04429v2}{arXiv:1909.04429v2}
(2026 revision of the 2019 preprint).

\bibitem{bbt} S.~Beckus, J.~Bellissard, and Y.~Thomas,
{\it Spectral regularity and defects for the Kohmoto model},
Annales Henri Poincar\'e \textbf{27} (2026), no.~6, 2117--2154.

\bibitem{bls}
J.~Bellissard, R.~Lima, and E.~Scoppola,
\textit{Localization in $\nu$-dimensional incommensurate
structures},
Comm. Math. Phys. \textbf{88} (1983), 465--477.

\bibitem{bs} J.~Bellissard and B.~Simon,
{\it Cantor spectrum for the almost Mathieu equation},
Journal of Functional Analysis \textbf{48} (1982), no.~3, 408-419.

\bibitem{bf} D.~S.~Borgnia and M.~Faust,
{\it The Dry Ten Martini Problem at Criticality},
preprint, arXiv:2112.06869v3
(2026 revision of the 2021 preprint).

\bibitem{cao}
H.~Cao,
\textit{Arithmetic localization for quasiperiodic long-range
operators with monotone potentials},
preprint, arXiv:2610.39448 (2026).

\bibitem{csz}
H.~Cao, Y.~Shi, and Z.~Zhang,
\textit{Localization for Lipschitz monotone quasi-periodic
Schr\"odinger operators on $\mathbb{Z}^d$ via Rellich
functions analysis},
Comm. Math. Phys. \textbf{406} (2025), Paper No.~102.

\bibitem{cl} C.~Cedzich and L.~Li,
{\it Twenty dry Martinis for the unitary almost-Mathieu operator},
Ergodic Theory and Dynamical Systems \textbf{46} (2026),
no.~4, 903--920.

\bibitem{ChoiElliottYui} M.-D.~Choi, G.~A.~Elliott, and N.~Yui,
{\it Gauss polynomials and the rotation algebra},
Inventiones Mathematicae \textbf{99} (1990), 225--246.

\bibitem{df} D.~Damanik and J.~Fillman,
{\it One-Dimensional Ergodic Schr\"odinger Operators:
I. General Theory},
Graduate Studies in Mathematics, vol.~221,
American Mathematical Society, Providence, RI, 2022.

\bibitem{DamanikKillip2005} D.~Damanik and R.~Killip,
{\it Ergodic potentials with a discontinuous sampling function are non-deterministic},
Mathematical Research Letters \textbf{12} (2005), no.~2, 187--192.

\bibitem{dl} 
D.~Damanik and D.~Lenz,
A condition of Boshernitzan and uniform convergence in the
multiplicative ergodic theorem,
\emph{Duke Math. J.} \textbf{133} (2006), 95--123.

\bibitem{TaoEtAl}
P.~B. Denton, S.~J. Parke, T. Tao, and X. Zhang,
{\it Eigenvectors from eigenvalues: A survey of a basic identity in linear algebra},
Bull. Amer. Math. Soc. (N.S.) \textbf{59} (2022), no.~1, 31--58.

\bibitem{dms} R.~del Rio, N.~Makarov, and B.~Simon,
{\it Operators with singular continuous spectrum:
II. Rank one operators},
Communications in Mathematical Physics \textbf{165} (1994),
no.~1, 59--67.

\bibitem{eli} L.~H.~Eliasson,
{\it Floquet solutions for the 1-dimensional quasi-periodic
Schr\"odinger equation},
Communications in Mathematical Physics \textbf{146} (1992), no.~3, 447--482.

\bibitem{fp}
A.~L.~Figotin and L.~A.~Pastur,
\textit{An exactly solvable model of a multidimensional
incommensurate structure},
Comm. Math. Phys. \textbf{95} (1984), 401--425.
\bibitem{ge}L.~Ge,
On the almost reducibility conjecture,
\emph{Geom. Funct. Anal.} \textbf{34} (2024), 32--59.
\bibitem{gjy24} L.~Ge, S.~Jitomirskaya, and J.~You,
{\it Kotani theory, Puig's argument, and stability of the Ten Martini Problem},
preprint, arXiv:2308.09321 (2023).

\bibitem{gj} L.~Ge and S.~Jitomirskaya,
{\it Intrinsic Symplectic Structure and Sharp Arithmetic Universality},
preprint, arXiv:2407.08866v2
(2026 revision of the 2024 preprint).

\bibitem{gjy} L.~Ge, S.~Jitomirskaya, and J.~You,
{\it The Robust Ten Martini Problem},
preprint, arXiv:2609.15812 (2026).

\bibitem{ge-dry} L.~Ge,
{\it The Dry Ten Martini Problem for subcritical Type I operators},
preprint, arXiv:2609.37438 (2026).

\bibitem{gjynew} L.~Ge, S.~Jitomirskaya, and J.~You,
{\it The Supercritical Ten Martini Problem},
preprint, arXiv:2611.13068 (2026).

\bibitem{gwx} L.~Ge, Y.~Wang, and J.~Xu,
{\it The Dry Ten Martini Problem for $C^2$ cosine-type
quasiperiodic Schr\"odinger operators},
preprint, arXiv:2503.06918 (2025).

\bibitem{gordon} A.~Ya.~Gordon,
{\it Pure point spectrum under 1-parameter perturbations
and instability of Anderson localization},
Communications in Mathematical Physics \textbf{164} (1994),
no.~3, 489--505.

\bibitem{gs} M.~Goldstein and W.~Schlag,
{\it On resonances and the formation of gaps in the spectrum
of quasi-periodic Schr\"odinger equations},
Annals of Mathematics (2) \textbf{173} (2011), no.~1, 337--475.

\bibitem{gfp}
D.~R.~Grempel, S.~Fishman, and R.~E.~Prange,
\textit{Localization in an incommensurate potential:
An exactly solvable model},
Phys. Rev. Lett. \textbf{49} (1982), 833--836.
\bibitem{hs}
B.~Helffer and J.~Sj\"ostrand,
Semi-classical analysis for Harper's equation. III.
Cantor structure of the spectrum,
\emph{M\'em. Soc. Math. France (N.S.)}
\textbf{39} (1989), 1--124.
\bibitem{hofstadter} D.~R.~Hofstadter,
{\it Energy levels and wave functions of Bloch electrons in rational
and irrational magnetic fields},
Physical Review B \textbf{14} (1976), no.~6, 2239--2249.

\bibitem{j} S.~Ya.~Jitomirskaya,
{\it Metal-insulator transition for the almost Mathieu operator},
Annals of Mathematics (2) \textbf{150} (1999), no.~3, 1159--1175.

\bibitem{JitomirskayaMavi14} S.~Jitomirskaya and R.~Mavi,
{\it Continuity of the measure of the spectrum for quasiperiodic Schr\"odinger operators with rough potentials},
Communications in Mathematical Physics \textbf{325} (2014), 585-601.

\bibitem{JitomirskayaMavi17} S.~Jitomirskaya and R.~Mavi,
{\it Dynamical bounds for quasiperiodic Schr\"odinger operators
with rough potentials},
International Mathematics Research Notices
\textbf{2017} (2017), no.~1, 96--120.

\bibitem{jk1} S.~Jitomirskaya and I.~Kachkovskiy, {\it All couplings localization for quasiperiodic operators with monotone potentials}, J.\ Eur.\ Math.\ Soc.\ \textbf{21} (2019), 777--795.

\bibitem{jk2} S.~Jitomirskaya and I.~Kachkovskiy, {\it Sharp arithmetic localization for quasiperiodic operators with monotone potentials}, preprint, arXiv:2407.00703.

\bibitem{jk3} S.~Jitomirskaya and I.~Kachkovskiy, {\it Sharp arithmetic delocalization for quasiperiodic operators with potentials of semi-bounded variation}, preprint. arXiv:2408.16935.

\bibitem{JitomirskayaKrasovsky2002} S.~Ya.~Jitomirskaya and I.~V.~Krasovsky,
{\it Continuity of the measure of the spectrum for discrete
	quasiperiodic operators},
Mathematical Research Letters \textbf{9} (2002), no.~4, 413--421.

\bibitem{jl} S.~Ya.~Jitomirskaya and Y.~Last,
{\it Anderson localization for the almost Mathieu equation, III.
Semi-uniform localization, continuity of gaps,
and measure of the spectrum},
Communications in Mathematical Physics \textbf{195} (1998),
no.~1, 1--14.

\bibitem{jls} S.~Jitomirskaya, N.~Levi, and T.~Stoiber,
{\it (Sharp arithmetic) localization beyond transversality},
preprint, arXiv:2610.17139 (2026).

\bibitem{jlmaryland}
S.~Jitomirskaya and W.~Liu,
\textit{Arithmetic spectral transitions for the Maryland model},
Comm. Pure Appl. Math. \textbf{70} (2017), no.~6, 1025--1051.

\bibitem{jliu} S.~Jitomirskaya and W.~Liu,
{\it Universal hierarchical structure of quasiperiodic eigenfunctions},
Annals of Mathematics (2) \textbf{187} (2018),
no.~3, 721--776.

\bibitem{jlt} S.~Jitomirskaya, W.~Liu, and S.~Tcheremchantsev,
{\it Lower bounds on concentration through Borel transforms and
quantitative singularity of spectral measures near the arithmetic
transition},
preprint, arXiv:2501.12153 (2025).


\bibitem{jymaryland}
S.~Jitomirskaya and F.~Yang,
\textit{Pure point spectrum for the Maryland model:
A constructive proof},
Ergodic Theory Dynam. Systems \textbf{41} (2021),
no.~1, 283--294.

\bibitem{jz} S.~Jitomirskaya and S.~Zhang,
{\it Quantitative continuity of singular continuous spectral
measures and arithmetic criteria for quasiperiodic
Schr\"odinger operators},
Journal of the European Mathematical Society \textbf{24}
(2022), no.~5, 1723--1767.

\bibitem{kach}
I.~Kachkovskiy,
\textit{Localization for quasiperiodic operators with unbounded
monotone potentials},
J. Funct. Anal. \textbf{277} (2019), no.~10, 3467--3490.

\bibitem{kpsseries}
I.~Kachkovskiy, L.~Parnovski, and R.~Shterenberg,
\textit{Convergence of perturbation series for unbounded
monotone quasiperiodic operators},
Adv. Math. \textbf{409} (2022), Paper No.~108647.

\bibitem{kpspert}
I.~Kachkovskiy, L.~Parnovski, and R.~Shterenberg,
\textit{Perturbative diagonalization and spectral gaps of
quasiperiodic operators on $\ell^2(\mathbb{Z}^d)$ with
monotone potentials},
Comm. Math. Phys. \textbf{406} (2025), Paper No.~125.

\bibitem{kkps}
I.~Kachkovskiy, S.~Krymski, L.~Parnovski, and R.~Shterenberg,
\textit{Perturbative diagonalization for Maryland-type
quasiperiodic operators with flat pieces},
J. Math. Phys. \textbf{62} (2021), no.~6, Paper No.~063509.

\bibitem{kps}
I.~Kachkovskiy, L.~Parnovski, and R.~Shterenberg,
{\it On gaps in the spectra of quasiperiodic Schr\"odinger
operators with discontinuous monotone potentials},
International Mathematics Research Notices
\textbf{2025}, no.~7, rnaf085 (2025).

\bibitem{kpsgeneric}
I.~Kachkovskiy, L.~Parnovski, and R.~Shterenberg,
\textit{Smooth perturbation theory and generic Cantor spectrum
for multi-dimensional monotone quasiperiodic operators},
preprint, arXiv:2610.00733 (2026).

\bibitem{kz}
J.~Kerdboon and X.~Zhu,
\textit{Anderson localization for Schr\"odinger operators
with monotone potentials over circle homeomorphisms},
J. Spectr. Theory \textbf{14} (2024), no.~4, 1623--1646.

\bibitem{last} Y.~Last,
{\it Zero measure spectrum for the almost Mathieu operator},
Communications in Mathematical Physics \textbf{164} (1994), no.~2, 421--432.

\bibitem{levi} N.~Levi,
{\it Universality of Packing Dimension Estimates for Spectral
Measures of Quasiperiodic Operators: Monotone Potentials},
preprint, arXiv:2502.03707 (2025).

\bibitem{lxz} X.~Li, D.~Xu, and Q.~Zhou,
{\it Monotonicity, global symplectification and the stability
of Dry Ten Martini Problem},
preprint, arXiv:2601.02222 (2026).

\bibitem{li-dry} X.~Li,
{\it Dry Ten Martini Problem in the Subcritical Type I Regime},
preprint, arXiv:2609.23467 (2026).

\bibitem{nuckolls2025}
K.~P.~Nuckolls, M.~G.~Scheer, D.~Wong, M.~Oh, R.~L.~Lee,
J.~Herzog-Arbeitman, K.~Watanabe, T.~Taniguchi, B.~Lian,
and A.~Yazdani,
\textit{Spectroscopy of the fractal Hofstadter energy spectrum},
Nature \textbf{639} (2025), 60--66.

\bibitem{Paige}
C.~C. Paige,
{\it The computation of eigenvalues and eigenvectors of very large sparse matrices},
Ph.D. thesis, University of London,
London University Institute of Computer Science, 1971.

\bibitem{peierls} R.~Peierls,
{\it Zur Theorie des Diamagnetismus von Leitungselektronen},
Zeitschrift f\"ur Physik \textbf{80} (1933), 763--791.

\bibitem{puig} J.~Puig,
{\it Cantor spectrum for the almost Mathieu operator},
Communications in Mathematical Physics \textbf{244} (2004), no.~2, 297--309.

\bibitem{simonmaryland}
B.~Simon,
\textit{Almost periodic Schr\"odinger operators.
IV. The Maryland model},
Ann. Physics \textbf{159} (1985), 157--183.

\bibitem{simon15} B.~Simon,
{\it Fifteen problems in mathematical physics},
in {\it Perspectives in Mathematics: Anniversary of Oberwolfach 1984},
W.~J\"ager, J.~Moser, and R.~Remmert (eds.),
Birkh\"auser, Basel, 1984, 423-454.

\bibitem{simonXXI} B.~Simon,
{\it Schr\"odinger operators in the twenty-first century},
in {\it Mathematical Physics 2000},
A.~Fokas, A.~Grigoryan, T.~Kibble, and B.~Zegarlinski (eds.),
Imperial College Press, London, 2000, 283-288.

\bibitem{sinai} Ya.~G.~Sinai,
{\it Anderson localization for one-dimensional difference
Schr\"odinger operator with quasiperiodic potential},
Journal of Statistical Physics \textbf{46} (1987), nos.~5--6, 861--909.

\end{thebibliography}
\end{document}